\documentclass[11pt]{amsart}
\usepackage{amsmath,amsthm, amscd, amssymb, amsfonts, mathrsfs,pb-diagram,pb-xy,mathtools}
\usepackage{multicol}
\usepackage[all]{xy}
\usepackage[inline]{enumitem}
\usepackage{mathrsfs}
\usepackage[dvips, dvipsnames, usenames]{color}
\usepackage{hyperref}
\usepackage[T2A,T1]{fontenc}
\usepackage{rotating}
\usepackage{multicol}
\usepackage{todonotes}

\usepackage{bbold}

\DeclareFontFamily{OT1}{pzc}{}
\DeclareFontShape{OT1}{pzc}{m}{it}{<-> s * [1.10] pzcmi7t}{}
\DeclareMathAlphabet{\mathpzc}{OT1}{pzc}{m}{it}

\allowdisplaybreaks

\newcommand{\trid}{\triangleright}
\newcommand{\fiz}{\triangleleft}

\newcommand{\uno}{\mathbf{1}}
\newcommand{\ghost}{\mathscr{G}}
\newcommand{\lstr}{\mathfrak L}

\newcommand{\ba}{ \mathbf{a}}
\newcommand{\bm}{ \mathbf{m}}
\newcommand{\bn}{ \mathbf{n}}
\newcommand{\bp}{\mathbf{p}}
\newcommand{\tb}{\mathbf{t}}

\newcommand{\sch}{\mbox{\usefont{T2A}{\rmdefault}{m}{n}\cyrshch}}

\newcommand{\sa}{\mathtt{r}}

\newcommand{\kut}{ \Bbbk^{\times}}
\newcommand{\Fp}{\mathbb F_{\hspace{-1pt} p}}

\makeatletter
\newcommand{\leqnomode}{\tagsleft@true}
\newcommand{\reqnomode}{\tagsleft@false}
\makeatother

\renewcommand{\_}[1]{_{\left( #1 \right)}}

\newcommand{\Dchaintwo}[3]{\xymatrix@C-4pt{\overset{#1}{\underset{\ }{\circ}}\ar
@{-}[r]^{#2}
& \overset{#3}{\underset{\ }{\circ}}}}

\newcommand{\Dchainfive}[9]{\xymatrix@C-6pt{\overset{#1}{\underset{\ }{\circ}}\ar
@{-}[r]^{#2}  & \overset{#3}{\underset{\ }{\circ}}\ar  @{-}[r]^{#4}  &
\overset{#5}{\underset{\ }{\circ}}
\ar  @{-}[r]^{#6}  & \overset{#7}{\underset{\ }{\circ}}\ar  @{-}[r]^{#8}  &
\overset{#9}{\underset{\ }{\circ}}}}

\newcommand{\mor}{\sim_{\text{Mor}}}

\newcommand{\Bg}{\mathfrak{B}}

\newcommand{\ab}{\mathfrak{ab}}
\newcommand{\ngo}{\mathfrak n}

\newcommand{\ad}{\operatorname{ad}}
\newcommand{\Alg}{\operatorname{Alg}}

\newcommand{\Der}{\operatorname{Der}}
\newcommand{\End}{\operatorname{End}}
\newcommand{\car}{\operatorname{char}}
\newcommand{\coh}{\operatorname{H}}

\newcommand{\gr}{\operatorname{gr}}
\newcommand{\Hom}{\operatorname{Hom}}
\newcommand{\Ext}{\operatorname{Ext}}
\newcommand{\id}{\operatorname{id}} 

\newcommand{\lcm}{\operatorname{lcm}}
\newcommand{\ord}{\operatorname{ord}}

\newcommand{\g}{\mathfrak g}
\newcommand{\Ggo}{\mathfrak{G}}
\newcommand{\lgo}{\mathfrak l}

\newcommand{\Lgo}{\mathfrak{L}}
\newcommand{\stt}{\mathtt{s}}
\newcommand{\ugo}{\mathfrak u}
\newcommand{\ago}{\mathfrak a}
\newcommand{\pgo}{\mathfrak p}
\newcommand{\rep}{\operatorname{rep}}
\newcommand{\corep}{\operatorname{corep}}
\newcommand{\diag}{\operatorname{diag}}

\newcommand{\ku}{ \Bbbk}
\newcommand{\fp}{{\mathbb F}_{\hspace{-2pt}p}}
\newcommand{\I}{\mathbb I}
\newcommand{\G}{\mathbb G}
\newcommand{\N}{\mathbb N}
\newcommand{\Z}{\mathbb Z}
\newcommand{\s}{\mathfrak{s}}

\newcommand{\toba}{\mathscr{B}}
\newcommand{\doublecpx}{\mathscr{E}}
\newcommand{\triple}{\mathscr T}

\newcommand{\bq}{\mathfrak{q}}

\newcommand{\Pc}{{\mathcal P}}
\newcommand{\Fc}{{\mathcal F}}

\newcommand{\Cc}{{\mathcal C}}

\newcommand{\Ss}{\mathcal{S}}

\newcommand{\yd}[1]{{}^{#1}_{#1}\mathcal{YD}}

\newcommand{\Ac}{\mathcal{A}}

\newcommand{\Rc}{\mathcal{R}}

\newcommand{\pf}{\begin{proof}}
\newcommand{\epf}{\end{proof}}

\newtheorem{theorem}{Theorem}[section]

\newtheorem{lemma}[theorem]{Lemma}

\newtheorem{cor}[theorem]{Corollary}

\newtheorem{prop}[theorem]{Proposition}

\theoremstyle{definition}
\newtheorem{definition}[theorem]{Definition}
\newtheorem{example}[theorem]{Example}

\newtheorem*{observation*}{Observation}
\newtheorem{question}[theorem]{Question}

\newtheorem*{claim*}{Claim}

\theoremstyle{remark}
\newtheorem{remark}[theorem]{Remark}

\newcommand{\ot}{\otimes}

\newcommand{\Huno}{\widetilde H}

\newcommand{\Vc}{\mathcal V}
\newcommand{\Vs}{\mathscr{V}}

\begin{document}

\noindent
\title[On the finite generation of the
cohomology of abelian extensions]
{On the finite generation of the
	cohomology of abelian extensions of  Hopf algebras}

\author[Andruskiewitsch and Natale]
{Nicol\'as Andruskiewitsch and Sonia Natale}

\address{FaMAF-CIEM (CONICET), Universidad Nacional de C\'ordoba,
Me\-dina A\-llen\-de s/n, Ciudad Universitaria, 5000 C\' ordoba,  Argentina.} 

\email{nicolas.andruskiewitsch@unc.edu.ar, natale@famaf.unc.edu.ar}

\thanks{  
The work of N. A. was  partially supported by Secyt (UNC), CONICET (PIP 11220200102916CO), FONCyT-ANPCyT (PICT-2019-03660) and 
 by the Shenzhen International Center of Mathematics, SUSTech. The research of S. N. was partially supported by Secyt (UNC)}

\begin{abstract}
A finite-dimensional Hopf algebra is called quasi-split if it is Morita equivalent
to a split abelian extension of Hopf algebras.
Combining results of Schauenburg and Negron, 
it is shown that  every quasi-split  finite-dimensional Hopf algebra 
satisfies the finite generation  cohomology conjecture of Etingof and Ostrik.
This is applied to  a 
family of pointed Hopf algebras in odd characteristic
introduced by Angiono, Heckenberger and the first author,
proving that they satisfy the aforementioned conjecture.
\end{abstract}

\maketitle

\setcounter{tocdepth}{1}

\section{Introduction}

\subsection{The problem}
Let $\ku$ be an algebraically closed field. 
We  say that an augmented $\ku$-algebra $A$ has finitely generated cohomology (fgc for short) 
when 

\smallbreak
\begin{enumerate}[leftmargin=*,label=\rm{(\alph*)}]
\item\label{item:fgca} the cohomology ring $\coh(A, \ku) = \oplus_{n\in \N_0} \Ext^n_A(\ku, \ku)$ is finitely generated, and 

\smallbreak
\item\label{item:fgcb} for any finitely generated $A$-module $M$, 
 $\coh(A,M) = \oplus_{n\in \N_0} \Ext^n_A(\ku, M)$ is a finitely generated $\coh(A, \ku)$-module. 
\end{enumerate}

\smallbreak
This definition was extended in \cite{etingof-ostrik03} as follows. A finite tensor category $\Cc$ 
(with unit object $\mathbb{1}$)
has finite generation of cohomology (or fgc) when

\smallbreak
\begin{enumerate}[leftmargin=*,label=\rm{(\alph*)}]
\item\label{item:fgc-tc-a} the $\ku$-algebra $\coh(\Cc) = \oplus_{n\in \N_0} \Ext^n(\mathbb{1}, \mathbb{1})$ 
is finitely generated, and 

\smallbreak
\item\label{item:fgc-tc-b} for any object $V$ in $\Cc$ the $\coh(\Cc) $-module  
$\coh(V) = \oplus_{n\in \N_0} \Ext^n(\mathbb{1}, V)$ is finitely generated.
\end{enumerate}

\medbreak
In  \cite{etingof-ostrik03}, Etingof and Ostrik, drawing on fundamental results of 
\cite{V,golod,evens,andersen-jantzen,
friedlander-parshall,friedlander-suslin,ginzburg-kumar}, conjectured that finite tensor categories have fgc,  in particular that finite-dimensional  Hopf algebras have fgc;
this was verified in many cases
\cite{bnpp,gordon,drupieski,drupieski2,MPSW,NWW,EOW,stefan-vay,friedlander-negron, negron,negron-plavnik,aapw}; see \cite[\S \, 1.1]{aapw} for background.

\smallbreak
\subsection{Morita equivalence}\label{subsec:morita}
Let $H$ and $U$ be finite-dimensional Hopf algebras.  
By  \cite[Theorem  3.4]{negron-plavnik}, see also \cite[Theorem 3.2.1]{aapw},  
if  the Drinfeld double $D(H)$ has fgc, then $H$ has fgc.
The following argument was intensively used in \cite{aapw}.

\medbreak
We  say  that  $H$ and  $U$ are
\emph{Morita equivalent}, denoted $H \mor U$, if the finite tensor categories $\rep H$ and $\rep U$ are Morita equivalent as in  \cite{muger,EGNO}. 
Equivalently, $H \mor U$ if there exists an equivalence of braided tensor categories 
between the Drinfeld centers $\mathcal Z(\rep H)$ and $\mathcal Z(\rep U)$  \cite{EGNO}.
In other words, $H \mor U$ iff  $D(H)$ and $D(U)$ are twist equivalent  quasitriangular Hopf algebras
(this does  not imply that $H$ and $U$ are Morita equivalent as algebras). 
Thus to prove that $H$  has fgc
it suffices to find   $U$ such that
\begin{align*}
\textrm{(i)} \ D(U) &\text{ has fgc}
& &\text{ and} &
\textrm{(ii)} \  &H \mor U.
\end{align*}

\subsection{Extensions}
In this paper we consider the question of fgc for a Hopf algebra $H$
 fitting into an exact sequence of finite-dimensional Hopf algebras
\begin{align}\label{eq:extensions-intro}
\ku \to K \rightarrow H \rightarrow L \to \ku.
\end{align}
If $H$ has fgc, then $K$ has fgc by \cite[Theorem 3.2.1]{aapw} while 
it is unclear whether we could infer that $L$ has fgc. Thus, it is natural to ask:

\begin{question}
Given an extension \eqref{eq:extensions-intro}  
such that $K$ and $L$ have fgc, does  $H$ also have fgc?
\end{question}

For instance, if $K$ is semisimple and $L$ has fgc, then $H$ has fgc, see \cite[Lemma 3.2.5]{aapw}.
The proof  uses a variation of the classical Hochschild-Serre spectral sequence
but it is not clear (to us) how to proceed when the semisimplicity assumption on $K$ is dropped.

\medbreak
The extension \eqref{eq:extensions-intro} is \emph{abelian} when $K$ is commutative 
and $L$ is cocommutative. In this case there are suitable
 actions of $L$ on $K$ and of $K^*$ on $L$, and a pair $(\sigma, \tau)$ of compatible cocycles  that determine the possible extensions $H$. 
 The suitable actions give rise to a double complex $\doublecpx$. 
 The abelian extension \eqref{eq:extensions-intro} is \emph{split} if 
 the pair $(\sigma, \tau)$ is trivial in the total complex associated to $\doublecpx$.
  Furthermore, by an argument due to G. I. Kac 
  there is a long exact sequence (the Kac exact sequence)
  in which the pair $(\sigma, \tau)$ is sent to a 3-cocycle.
  See \cite{masuoka-survey} for details.
 
We shall say that an abelian extension is \emph{quasi-split} 
if it is Morita equivalent to the split extension (with respect to the same suitable actions);
see Definition \ref{def:quasi-split}.

Tautologically, split abelian extensions are quasi-split; also, abelian extensions for which the Kac 3-cocycle is trivial are quasi-split (see Corollary \ref{cor:quasi-split-phitrivial}). The starting 
point of the paper is the following result.

\begin{theorem}\label{thm:fingencoh-extensions}
If \eqref{eq:extensions-intro} is a quasi-split abelian extension, then $H$ and $D(H)$ have fgc.
\end{theorem}

See Theorem \ref{thm:singer-pair} for a precise formulation.
As a consequence, the dual, any twist and any cocycle deformation
of $H$ have fgc. In characteristic 0, Theorem \ref{thm:fingencoh-extensions}
is trivial, as a finite-dimensional abelian extension is semisimple.

\medbreak
Theorem \ref{thm:fingencoh-extensions} follows combining 
the next two results due to Negron and Schauenburg respectively.

\begin{theorem}\label{thm:negron} \cite{negron}
If $U$ is a cocommutative Hopf algebra, $\dim U < \infty$, then  $D(U)$  has fgc.
\end{theorem}

\medbreak
Notice that Theorem \ref{thm:negron} is stated in  \cite{negron} in the language of finite group schemes.

\medbreak
\begin{theorem}\label{thm:abelian-extensions} \cite{schauenburg}
Let \eqref{eq:extensions-intro} be a split abelian extension 
of  finite-dimensional Hopf algebras. Then there exists a cocommutative Hopf algebra
$U$ such that $H \mor U$.
\end{theorem}
Observe that the Drinfeld double of a finite-dimensional  
cocommutative Hopf algebra is a split abelian extension.

\medbreak
\subsection{Pointed Hopf algebras} \label{sub:applications} 
In the paper \cite{aah-triang} (assuming $\car \ku =0$)
Nichols algebras over abelian groups 
with finite Gelfand-Kirillov dimension of a certain kind were classified. 
It was observed later in  \cite{aah-oddchar}
that many of the Nichols algebras over the analogous
 braided vector spaces in odd characteristic have finite dimension,
 hence give rise to finite dimensional pointed Hopf algebras by bosonization
 with group algebras of finite abelian groups. 
Theorem \ref{thm:restricted-gen-ext}, the main result of this paper,  
shows that  these  pointed Hopf algebras  have fgc when the groups are well
chosen. 
Concretely, in the notation of Section \ref{sec:hopf}, we consider a family of
braided vector spaces $\Vs(\bq, \ba)$, we fix a suitable 
finite group $\Gamma$ such that the Nichols algebra $\toba(\Vs(\uno, \ba))$
can be realized as a Hopf algebra in the category $\yd{\ku \Gamma}$
of Yetter-Drinfeld modules over $\Gamma$ and proceed in two stages:

\begin{enumerate}[leftmargin=*,label=\rm{(\roman*)}]
\item we show that the bosonization $H = \toba(\Vs(\uno, \ba)) \# \ku \Gamma$
fits into a \emph{split} abelian exact sequence, hence $D(H)$ and $H$ have fgc by Theorem \ref{thm:fingencoh-extensions};

\medbreak
\item for a general $\bq$, we present $\toba(\Vs(\bq, \ba)) \# \ku \Gamma$
as a cocycle deformation of $H$, hence $\toba(\Vs(\bq, \ba)) \# \ku \Gamma$ 
has fgc too.
\end{enumerate}
 We observe:
 \begin{itemize}[leftmargin=*]\renewcommand{\labelitemi}{$\circ$}
 \item Not  all  Nichols algebras in \cite{aah-oddchar} belong  to the family
 treated here;  most of the remaining ones
 arise from abelian extensions of Hopf superalgebras and will be dealt with elsewhere.
 
 \medbreak 
 \item Many realizations of $\toba(\Vs(\uno, \ba))$ fit into  abelian 
 exact sequences which are not split.
 	
 \medbreak 
 \item Many realizations of $\toba(\Vs(\bq, \ba))$ fit into  exact sequences
 	which are not abelian.
 \end{itemize}

\subsection{Organization}\label{subsection:organization}
Abelian extensions are discussed in Section \ref{sec:extensions}.
 Section \ref{sect:quasi-split} is devoted to 
 Theorem \ref{thm:abelian-extensions}.
 The analysis of the Nichols algebras $\toba(\Vs(\bq, \ba))$ 
 is delicate as it involves a number of distinct combinatorial
 features developped in \cite{aah-triang}. For clarity, we
 first deal  with the two simplest examples, namely the restricted Jordan 
 plane and the first Laestrygonian algebra, in Sections \ref{sec:Jordan-plane}
 and \ref{section:lstr-11disc} respectively.
 In Section \ref{sec:hopf}, we work out the strategy outlined above for the 
 general Nichols algebras $\toba(\Vs(\bq, \ba))$.

\subsection{Conventions}\label{subsection:conventions}
For $\ell < \theta \in\N_0$, we set $\I_{\ell, \theta}=\{\ell, \ell +1,\dots,\theta\}$, $\I_\theta 
= \I_{1, \theta}$. 
Let $\G_N$ be the group of roots of unity of order $N$ in $\ku$ and $\G_N'$ the subset of primitive roots of order $N$;
$\G_{\infty} = \bigcup_{N\in \N} \G_N$ and $\G'_{\infty} = \G_{\infty} - \{1\}$.
If $L \in \N$ and $q\in \ku^{\times}$, then $(L)_q \coloneqq  \sum_{j=0}^{L-1}q^{j}$.

All  vector spaces, algebras and tensor products  are over $\ku$. 
We use $V^*$ to denote the linear dual of a vector space $V$, $V^* = \Hom_k(V,k)$. 

By abuse of notation, $\langle a_i: i\in I\rangle$ denotes either the group, the subgroup or the vector subspace generated by all $a_i$ for $i$ in an indexing set $I$,
the meaning being clear from the context. Instead, the subalgebra generated by all $a_i$ for $i \in I$ is denoted by $\ku \langle a_i: i\in I\rangle$.

The notation for Hopf algebras is standard: $\Delta$, $\varepsilon$, $\Ss$
denote the comultiplication, the counit, the antipode (always assumed bijective), 
respectively.
We use the Sweedler notation for the comultiplication and the coactions.
Our reference for the theory of Hopf algebras is \cite{radford-book}.
Generalities on Nichols algebras can be found in  \cite{andrus-leyva}.

\section{Extensions of Hopf algebras}\label{sec:extensions}
This Section contains a crash exposition of extensions of Hopf algebras.

\subsection{Exact sequences}
Recall that a sequence of morphisms of Hopf algebras
\begin{align}\label{eq:exact}
\ku \to A\overset{\iota}\to C \overset{\pi}\to B \to \ku
\end{align}
is exact \cite{andrus-devoto, schneider,hofstetter} if the following conditions holds:
\begin{multicols}{2}
\begin{enumerate}[leftmargin=*,label=\rm{(\roman*)}]
\item\label{suc-exacta-1} $\iota$ is injective.
\item\label{suc-exacta-2} $\pi$ is surjective.
\item\label{suc-exacta-3} $\ker\pi = C\iota(A)^+$.
\item\label{suc-exacta-4} $\iota(A) = C^{\operatorname{co} \pi}$.
\end{enumerate}
\end{multicols}

The exact sequence \eqref{eq:exact} is \emph{abelian} if $A$ is commutative and $B$ is cocommutative.
We shall also refer to $C$ in \eqref{eq:exact} as an (abelian when it corresponds)
extension of $B$ by $A$.

\medbreak 
The exact sequence \eqref{eq:exact} is \emph{cleft} 
(see \cite[Definition 3.2.13]{andrus-devoto}) if there exist 
\begin{enumerate}[leftmargin=*,label=\rm{(\roman*)}]
\item a unit preserving right $B$-colinear section $\mathpzc{s}: B \to C$ of the  projection $\pi$  and

\smallbreak
\item  a counit preserving left $A$-linear retraction $\mathpzc{r}: C \to A$ of the inclusion $\iota$,
\end{enumerate}
 both $\mathpzc{s}$ and $\mathpzc{r}$ being invertible with respect to the convolution product,
such that  the following equivalent conditions hold, for all $c \in C$:

\medbreak 
\begin{enumerate}[leftmargin=*,label=\rm{(\alph*)}]
\item $\mathpzc{s}^{-1}(\pi(c)) = \mathcal S(c_{(1)}) \mathpzc{r} (c_{(2)})$.

\smallbreak
\item $\mathpzc{s}(\pi(c)) = \mathpzc{r}^{-1}(c_{(1)}) c_{(2)}$.

\smallbreak
\item $\mathpzc{r}^{-1} (c) = \mathpzc{s}(\pi(c_{(1)})) \mathcal S(c_{(2)})$.

\smallbreak
\item $\mathpzc{r}(c) = c_{(1)} \mathpzc{s}^{-1} (\pi(c_{(2)}))$.

\smallbreak
\item $\mathpzc{r} \mathpzc{s} = \epsilon_B1_A$.
\end{enumerate}

The maps $\mathpzc{s}$ and $\mathpzc{r}$ are called \emph{cleaving maps} of \eqref{eq:exact}. 

\begin{remark} \cite{schneider}
If  $C$ is finite dimensional, then  the extension \eqref{eq:exact} is cleft. 
\end{remark}

\begin{definition}\cite[Definition 6.5.2]{schauenburg}
The exact sequence \eqref{eq:exact} is called \emph{split} if there exist  cleaving maps $\mathpzc{s}$ and $\mathpzc{r}$ as above, called \emph{splittings} of \eqref{eq:exact}, such that $\mathpzc{s}$ is an algebra map and $\mathpzc{r}$ is a coalgebra map.
\end{definition}

\medbreak 
Cleft extensions can be described using suitable linear maps.
A \emph{bicrossed product datum} of Hopf algebras is a collection 
$(A, B, \rightharpoonup, \rho, \sigma, \tau)$ 
where $A$ and $B$ are Hopf algebras, 
\begin{align*}
\rightharpoonup: B \otimes A &\to A, &
\rho: B &\to B \otimes A,& 
\sigma: B \otimes B &\to A &&\text{and} & \tau: B &\to A \otimes A 
\end{align*}
are maps, called, respectively, a weak action, a weak coaction, a cocycle and a dual cocycle, obeying the conditions in \cite[Theorem 2.20]{andrus-devoto}. 

\medbreak 
A bicrossed product datum $(A, B, \rightharpoonup, \rho, \sigma, \tau)$  gives rise to a Hopf algebra $A \#{}_{\sigma}^{\tau}B$, called a \emph{bicrossed product}: the underlying vector space is $A \otimes B$, while the multiplication, comultiplication and antipode are determined by the  formulas
\begin{align*}
(k \# h)(t \# g) &= k\left(h_{(1)} \rightharpoonup t\right) \sigma(h_{(2)}, g_{(1)}) \# h_{(3)}  g_{(2)},
\\
\Delta(k \# h) &= k_{(1)} \tau^{(1)}(h_{(1)}) \# \rho(h_{(2)})_i \otimes k_{(2)} \tau^{(2)}(h_{(1)})\rho(h_{(2)})^i \#  h_{(3)},
\\
\mathcal S(a \# b) & = ((\mathcal S((\rho(b)_i)_{(2)}) \rightharpoonup \mathcal S(\rho(b)^i) \# \mathcal S((\rho(b)_i)_{(1)})) \mathcal S(a) \# 1.
\end{align*}

The natural maps $\iota: A \to A \#{}_{\sigma}^{\tau}B$ and $\pi: A \#{}_{\sigma}^{\tau}B \to B$ fit into a cleft exact sequence 
$$\ku \to A\overset{\iota}\to  A \#{}_{\sigma}^{\tau}B  \overset{\pi}\to B \to \ku,$$ 
of Hopf algebras, with cleaving maps \begin{align*}
\mathpzc{s}: B &\to A \#{}_{\sigma}^{\tau}B, & \mathpzc{s}(b) &= 1 \# b, &
\mathpzc{r}: A \#{}_{\sigma}^{\tau}B &\to A, & \mathpzc{r}(a \# b) &= \epsilon(b) a,&  a &\in A, \, b \in B. 
\end{align*}
\medbreak 
Conversely, given a cleft exact sequence \eqref{eq:exact}, there exists a bicrossed product datum $(A, B, \rightharpoonup, \rho, \sigma, \tau)$ such that $C \cong A \#{}_{\sigma}^{\tau}B$; it arises from the cleaving maps $\mathpzc{s}$ and $\mathpzc{r}$ as follows:
\begin{align*}
b \rightharpoonup a &= \mathpzc{s}(b_{(1)}) a \mathpzc{s}^{-1}(b_{(2)}),
\quad 
\sigma(b \# b') = \mathpzc{s}(b_{(1)}) \mathpzc{s}(b'_{(1)}) \mathpzc{s}^{-1}(b_{(2)} b'_{(2)}),
\\
\rho(\pi(c)) & = \pi(c_{(2)}) \otimes \mathpzc{r}^{-1}(c_{(1)}) \mathpzc{r}(c_{(3)}),
\quad 
\tau (\pi(c))  =  \Delta(\mathpzc{r}^{-1}(c_{(1)})) \; 
\mathpzc{r}(c_{(2)}) \otimes \mathpzc{r}(c_{(3)}),
\end{align*}
for all $a\in A$, $b \in B$, $c \in C$. See \cite[Theorem 3.2.14]{andrus-devoto}.

\medbreak 
A special instance of a bicrossed product datum $(A, B, \rightharpoonup, \rho, \sigma, \tau)$ 
occurs where the cocycle $\sigma$ and the dual cocycle $\tau$ are trivial maps; 
in this case, we omit the mention of $\sigma$ and $\tau$ and call 
the bicrossed product datum $(A, B, \rightharpoonup, \rho)$ a \emph{bismash datum}. 
Notice that in a bismash datum,  $A$ is a left $B$-module algebra and $B$ is a right $A$-comodule algebra with  action $\rightharpoonup$ and coaction $\rho$.

\medbreak
Given a bismash datum $(A, B, \rightharpoonup, \rho)$ ,
the   associated bicrossed product is denoted $A \# B$. The canonical cleaving maps imply that the associated exact sequence 
\begin{equation}\label{split-ab} \ku \to A\overset{\iota}\to A\#B \overset{\pi}\to B \to \ku\end{equation} 
is  split. Conversely, if \eqref{eq:exact} is a split exact sequence, then the corresponding bicrossed product datum $(A, B, \rightharpoonup, \rho, \sigma, \tau)$ is in fact a bismash datum $(A, B, \rightharpoonup, \rho)$ such that $C \cong A\# B$.

\medbreak
A bismash datum $(A, B, \rightharpoonup, \rho)$ is called a 
\emph{Singer pair} if $A$ is commutative and $B$ is cocommutative. In this case 
\eqref{split-ab} is a split abelian extension of Hopf algebras.
Observe that every abelian cleft exact sequence \eqref{eq:exact} gives rise to a Singer pair through the actions $\rightharpoonup$, $\rho$.

\subsection{Matched pairs}\label{sec:matched-pairs}
We now present a way to produce extensions due to G. I. Kac for abelian extensions
and to S. Majid in general. We start by the definition, see  \cite[\S 7.2]{majid-book}. 
An \emph{exact factorization} of a Hopf algebra $S$ consists of a pair $(G,L)$ of Hopf subalgebras of $S$
such that the restriction of the multiplication map 
\begin{align*}
\operatorname{mult}: G \otimes L \to S
\end{align*}
is a linear isomorphism. Exact factorizations  are classified through the following notion.

\begin{definition} \cite[7.2.1]{majid-book}
A \emph{matched pair} of Hopf algebras is a collection $(G, L, \trid, \fiz)$ where $L$ and $G$ Hopf algebras,
$G$ is a left $L$-module coalgebra with action $\trid$, 
$L$ is a right $G$-module coalgebra with action $\fiz$  such that for all $\ell,m \in L$, $x, y \in G$:
\begin{align}\label{eq:mp-hopf1}
(\ell m) \fiz x &= \left(\ell \fiz (m_{(1)} \trid x_{(1)})\right)(m_{(2)} \fiz x_{(2)}),
\\ \label{eq:mp-hopf2}
\ell \trid (x y ) &=\left(\ell_{(1)} \trid x_{(1)}\right)\left((\ell_{(2)} \fiz x_{(2)}) \trid y \right),
\\ \label{eq:mp-hopf3}
\ell_{(1)} \fiz  x_{(1)} \otimes \ell_{(2)} \trid x_{(2)} &= \ell_{(2)} \fiz x_{(2)} \otimes \ell_{(1)} \trid x_{(1)}.
\end{align} 
\end{definition}

\begin{prop} \label{prop:matchedpair-hopf-algebras-exact-factorization}
\begin{enumerate}[leftmargin=*,label=\rm{(\roman*)}]
\item \cite[7.2.2]{majid-book}
Given a matched pair  $(G, L, \trid, \fiz)$,  the coalgebra $G \otimes L$ with the multiplication
\begin{align*}
(x \otimes \ell)(y \otimes m) &= x\left(\ell_{(1)} \trid y_{(1)}\right) \otimes\left(\ell_{(2)} \fiz y_{(2)}\right) m,
& \ell,m &\in L, \ x, y \in G,
\end{align*}
\end{enumerate}
\noindent is a Hopf algebra denoted $G \bowtie L$. The natural inclusions $G \hookrightarrow G \bowtie L \hookleftarrow L$
form an exact factorization of $G \bowtie L$.

\medbreak
\begin{enumerate}[resume, leftmargin=*,label=\rm{(\roman*)}]
\item \cite[7.2.3]{majid-book} If $(G, L)$ is an exact factorization of a Hopf algebra $S$, then
there are actions $\trid$ and $\fiz$ such that $(G, L, \trid, \fiz)$ is a matched pair and 
\begin{align*}
S &\simeq G \bowtie L. 
\end{align*}\end{enumerate}
\end{prop}

The Hopf algebra $G \bowtie L$ is called the \emph{double crossproduct} associated to the actions $\fiz, \trid$. Matched pairs of Hopf algebras give rise to split exact sequences 
in the following way.

\begin{prop}\label{prop:matched-pairs-hopf-classic} 
\cite[7.2.3]{majid-book} If $(G, L, \trid, \fiz)$ is a matched pair of a Hopf algebras 
such that $\dim L < \infty$, 
then $(L^{*}, G, \rightharpoonup, \rho)$ is a bismash datum, where $\rightharpoonup$ 
and $\rho$ are obtained  by dualization; and vice versa.
\end{prop}

In conclusion, finite-dimensional split abelian extensions are determined by exact factorizations
of \emph{cocommutative} Hopf algebras  (that can be thought of as finite group schemes).
We illustrate these notions with some examples, see Examples \ref{exa:mp-finite-gps} and \ref{exa:mp-restricted-lie} for a discussion of the
fgc property.

\medbreak
\subsubsection*{\it Group algebras}
The exact factorizations of a group algebra $\ku \varSigma$ are in bijective correspondence 
with the exact factorizations of the group $\varSigma$, \cite{kac,mackey,Takeuchi-matched}. 
Also, the matched pairs of Hopf algebras $(\ku \varGamma, \ku \varLambda, \trid, \fiz)$ 
are the linearizations of the matched pairs 
of groups $(\varGamma, \varLambda, \trid, \fiz)$.

\medbreak
\subsubsection*{\it Lie algebras}
The exact factorizations of an enveloping algebra $U(\s)$ are in bijective correspondence 
with those of the Lie algebra $\s$;  matched pairs 
of Hopf algebras $(U(\g), U(\lgo), \trid, \fiz)$ are in bijective correspondence with matched pairs 
of Lie algebras $(\g, \lgo, \trid, \fiz)$. See \cite[\S 8.3]{majid-book}.
 More precisely,
 \begin{itemize}  [leftmargin=*]\renewcommand{\labelitemi}{$\circ$}
\item an exact factorization of a Lie algebra $\s$ consists of a pair $(\g, \lgo)$ of Lie
subalgebras such that $\s = \g \oplus \lgo$ (as vector spaces);

\item a \emph{matched pair} of Lie algebras is a collection $(\g, \lgo, \trid, \fiz)$ 
where $\g$ and $\lgo$ are Lie algebras,
$\trid$ and $\fiz$ are left and right actions 
$\xymatrix{\lgo & \lgo \times \g \ar  @{->}[r]^{\trid }\ar  @{->}[l]_{\fiz } & \g }$
satisfying \eqref{eq:mp-liealg2} and \eqref{eq:mp-liealg1} below.
 \end{itemize}
 
 \smallbreak
Given such a matched pair, $\g \bowtie \lgo \coloneqq  \g \oplus \lgo$ 
with the multiplication given by \eqref{eq:mp-liealg-bracket} is a Lie algebra. 
Up to identifications,  $(\g, \lgo)$ is an exact factorization of $\g \bowtie \lgo$. 
Also, \eqref{eq:mp-liealg-bracket} is equivalent to $\g$ and $\lgo$ being Lie subalgebras and 
\begin{align}\label{eq:mp-liealg-bracket2}
\left[\ell, y\right] &= \ell \trid y +  \ell \fiz y, & \ell &\in \lgo, \ y \in \g.
\end{align}
Conversely if $(\g, \lgo)$ is an exact factorization of a Lie algebra $\s$, then  \eqref{eq:mp-liealg-bracket2} defines the actions $\trid$ and $\fiz$,  
$(\g, \lgo, \trid, \fiz)$  is a matched pair and $\s \simeq \g \bowtie \lgo$.

\medbreak
\subsubsection*{\it Restricted Lie algebras}
In this example $\car \ku = p > 0$. 
Let $\s$ be a restricted Lie algebra 
with  $p$-operation $s\mapsto s^{[p]}$, see \cite[(65), p. 187]{jacobson}.
Concretely, given $i \in \I_{p-1}$, recall that
$\stt_{i}: \s \times \s \to \s$ is the homogeneous polynomial of degree $p$
defined    on a pair
$(s,t)$ as the coefficient of $X^{i-1}$ in  $\ad(Xs+t)^{p-1}(s)$, where $X$ is a formal variable. 
Then
\begin{align}\label{eq:restricted0}
(s+t)^{p} &= s^{p}+t^{p}+\sum _{i=1}^{p-1}\dfrac {\stt_{i}(s,t)}{i}, & & s, t\in \s.
\end{align}
By definition the $p$-operation satisfies for all $k\in \ku$, $s,t\in \s$ the identities
\begin{align}\label{eq:restricted1}
(ks)^{[p]} &= k^{p}s^{[p]}, 
\\\label{eq:restricted2}
\ad(s^{[p]}) &= \ad(s)^{p},
\\ \label{eq:restricted3}
(s+t)^{[p]} &= s^{[p]}+t^{[p]}+\sum _{i=1}^{p-1}\dfrac {\stt_{i}(s,t)}{i}.
\end{align}

The following definitions are natural.
 \begin{itemize}  [leftmargin=*]\renewcommand{\labelitemi}{$\circ$}
\item  An \emph{exact factorization} of  $\s$ is a pair $(\g, \lgo)$ of restricted Lie
subalgebras such that $\s = \g \oplus \lgo$ (as vector spaces). 

\medbreak
\item A \emph{matched pair} of restricted Lie algebras is a collection $(\g, \lgo, \trid, \fiz)$ 
where $\g$ and $\lgo$ are restricted Lie algebras,
$\trid$ and $\fiz$ are left and right $p$-actions 
$\xymatrix{\lgo & \lgo \times \g \ar  @{->}[r]^{\trid }\ar  @{->}[l]_{\fiz } & \g }$
satisfying for all $\ell,m  \in \lgo$, $x, y \in \g$ the identities
\begin{align}
\label{eq:mp-liealg2}
[\ell, m]\fiz x &= [\ell \fiz x, m]+[\ell, m \fiz x]+\ell \fiz(m  \trid  x) - m \fiz(\ell \trid x),
\\
\label{eq:mp-liealg1}
\ell \trid [x, y] &= \left[\ell \trid  x, y\right]+[x, \ell \trid y]+(\ell \fiz x)  \trid  y-(\ell \fiz y)  \trid  x.
\\
\label{eq:mp-restricted-liealg1}
\ell^{[p]} \fiz y  &= 
\sum_{1 \le i \le p-1} (\ad \ell)^{i} \left(\ell \fiz(\ell^{p-i} \trid y)\right), 
\\ \label{eq:mp-restricted-liealg2}
\ell \trid y^{[p]}  &= 
\sum_{1 \le i \le p-1} (-1)^{p-i}(\ad y)^{i} \left((\ell \fiz y^{p-1-i}) \trid y\right).
\end{align}
\end{itemize}

\begin{lemma} \label{lemma:mp-restricted}
Let $\g$ and $\lgo$ be restricted Lie algebras.
\begin{enumerate}[leftmargin=*,label=\rm{(\roman*)}]
\item \label{item:mp-restricted1}
If $(\g, \lgo, \trid, \fiz)$ is a matched pair of restricted Lie algebras,
then $\g \bowtie \lgo \coloneqq  \g \oplus \lgo$ is a restricted Lie algebra   with the multiplication
\begin{align}\label{eq:mp-liealg-bracket}
\left[(x , \ell), (y , m)\right] &= \left([x,y] + \ell \trid y - m \trid x,  
[\ell, m] + \ell \fiz y - m \fiz x \right)
\end{align}
$ \ell, m \in \lgo$, $x, y \in \g$; and with $p$-operation extending those
of $\g$ and $\lgo$ and 
\begin{align}\label{eq:mp-p-operation}
(y+\ell)^{[p]} &= y^{[p]}+\ell^{[p]}+\sum _{i=1}^{p-1}\dfrac {\stt_{i}(y,\ell)}{i}, &
y \in \g,\  &\ell \in \lgo,
\end{align}
where $\stt_{i}: \g \bowtie \lgo \times \g \bowtie \lgo \to \g \bowtie \lgo$ is defined from the Lie bracket \eqref{eq:mp-liealg-bracket}.

\item \label{item:mp-restricted2}
Let $(\g, \lgo)$ be an exact factorization of a restricted Lie algebra $\s$.
Then $(\g, \lgo, \trid, \fiz)$ with the actions given by \eqref{eq:mp-liealg-bracket2} is a matched pair of restricted Lie algebras and $\s \simeq \g \bowtie \lgo$ as restricted Lie algebras.
\end{enumerate}
\end{lemma}

\pf \ref{item:mp-restricted1}:  Since  $\trid, \fiz$  are  actions that satisfy  \eqref{eq:mp-liealg2} and  \eqref{eq:mp-liealg1}, 
$\g \bowtie \lgo$ with the bracket \eqref{eq:mp-liealg-bracket} is a Lie algebra \cite[\S 8.3]{majid-book}.
Thus the maps $\stt_i$ are defined and
we just have to check that \eqref{eq:mp-p-operation} gives a $p$-operation. 
\eqref{eq:restricted1} holds because $\stt_i$ is homogeneous of degree $p$.
We first verify \eqref{eq:restricted2} for $s = y \in \g$. Since both sides are linear operators 
and the restriction to $\g$ is a $p$-operation, it is enough to see that 
$\ad(y^{[p]}) (\ell) \overset{\star}{=} \ad(y)^{p}(\ell)$ for $\ell \in \lgo$.
Arguing by induction from \eqref{eq:mp-liealg-bracket2} we prove that for every $N \in \N$ and $\ell \in \lgo$
\begin{align*}
(\ad \ell)^N (y) &= \ell^N  \trid y + \sum_{0 \le i \le N-1} (\ad \ell)^i \left(y \fiz (\ell^{N-1-i}  \trid y ) \right).
\end{align*}
Taking $N=p$ and again by \eqref{eq:mp-liealg-bracket2}, we conclude that $\star$ holds iff \eqref{eq:mp-restricted-liealg1} is true. Similarly we prove that for every $N \in \N$, $y \in \g$ and $\ell \in \lgo$
\begin{align*}
(\ad y)^N (\ell) &=  \sum_{0 \le i \le N-1} (-1)^{N+ i} (\ad y)^i \left( (\ell \fiz y^{N-1-i})  \trid y ) \right)
+(-1)^N \ell  \fiz y^N.
\end{align*}
(Here \eqref{eq:mp-liealg-bracket2} says that $(\ad y) (\ell) = -\ell \trid y - \ell \fiz y$).
Taking $N=p$, we conclude that 
$\ad(\ell^{[p]}) (y) = \ad(\ell)^{p}(y)$ holds iff \eqref{eq:mp-restricted-liealg2} is true. 
Finally 
\begin{align*}
\ad(y + \ell)^{p} &= \ad (y)^{p}+ \ad(\ell)^{p} + \sum _{i=1}^{p-1}\dfrac {\ad \stt_{i}(y,\ell)}{i} 
\\ &= \ad (y^{[p]}) + \ad(\ell^{[p]}) + \sum _{i=1}^{p-1}\dfrac {\ad \stt_{i}(y,\ell)}{i}
= \ad ((y + \ell)^{[p]}),
\end{align*}
where the first equality is by \eqref{eq:restricted0}, the second because these are $p$-operations on the subalgebras and the third is by definition.
Thus \eqref{eq:restricted2} holds.

We proceed with \eqref{eq:restricted3}. Let $s = x + \ell$ and $t = y + m$, where $x,y \in \g$, $\ell, m \in \lgo$. Then
\begin{align*}
(s+t)^{[p]} &= \left((x+y) + (\ell + m)\right)^{[p]} 
\\
&= (x+y)^{[p]}+(\ell + m)^{[p]}+\sum _{i=1}^{p-1}\dfrac {\stt_{i}((x+y),(\ell + m))}{i}, 
\\ &= x^{[p]} + y^{[p]}+\sum _{i=1}^{p-1}\dfrac {\stt_{i}(x,y)}{i} 
+ \ell^{[p]}+ m^{[p]}+\sum _{i=1}^{p-1}\dfrac {\stt_{i}(\ell, m)}{i}
\\ & \qquad +\sum _{i=1}^{p-1}\dfrac {\stt_{i}((x+y),(\ell + m))}{i}.
\end{align*}
On the other hand,
\begin{align*}
s^{[p]}+t^{[p]} + \sum _{i=1}^{p-1}\dfrac {\stt_{i}(s,t)}{i}  &=
x^{[p]} + \ell^{[p]} + \sum _{i=1}^{p-1}\dfrac {\stt_{i}(x, \ell)}{i}
\\ &+
y^{[p]}+m^{[p]} + \sum _{i=1}^{p-1}\dfrac {\stt_{i}(y,m)}{i}
+ \sum _{i=1}^{p-1}\dfrac {\stt_{i}(s,t)}{i}
\end{align*}
To show that these two expressions are equal just apply \eqref{eq:restricted0}
to both sides of the equality
$\left((x+y) + (\ell + m) \right)^p = \left((x+ \ell) + (y + m) \right)^p$.
\epf

Let $\ugo(\s)$ be the restricted enveloping algebra of the restricted Lie algebra $\s$; it has a PBW-basis
with the powers of the generators truncated at $p$.

\begin{lemma} \label{lemma:exfact-restricted}
Let $\s$ be a restricted Lie algebra. The following are equivalent: 
\begin{enumerate}[leftmargin=*,label=\rm{(\roman*)}]
\item \label{item:pexfact-restricted1}
Exact factorizations of the Hopf algebra $\ugo(\s)$.

\item \label{item:pexfact-restricted2}
Exact factorizations of the restricted Lie algebra $\s$.
\end{enumerate}
\end{lemma}

\pf \ref{item:pexfact-restricted1} $\Rightarrow$ \ref{item:pexfact-restricted2}.
If $(\Ggo, \Lgo)$ is an exact factorization of $\ugo(\s)$, then take 
$\g = \Pc(\Ggo)$, $\lgo = \Pc(\Lgo)$. As a consequence of
\cite[Proposition 13.2.3]{sweedler}, $\Ggo \simeq \ugo(\g)$,
$\Lgo  \simeq \ugo(\lgo)$, and by the PBW theorem, 
$(\g, \lgo)$ is an exact factorization of $\s$.

\ref{item:pexfact-restricted1} $\Leftarrow$ \ref{item:pexfact-restricted2}: 
If $(\g, \lgo)$ is an exact factorization of $\s$, then the multiplication
$\ugo(\g) \otimes \ugo(\lgo)\to \ugo(\s)$ is a linear isomorphism--apply the PBW-theorem 
to the union of bases of $\g$ and $\lgo$.
\epf

Here are some concrete examples of factorizations of restricted Lie algebras.

\subsubsection*{\it Restricted Lie bialgebras}
Let us say that a  finite-dimensional Lie bialgebra $\mathfrak b$ is restricted if and only if its Manin triple 
$(\pgo, \mathfrak b, \mathfrak b^*)$ is restricted, meaning that $\pgo$ is restricted and 
$\mathfrak b$, $\mathfrak b^*$ are restricted subalgebras. In particular $(\mathfrak b, \mathfrak b^*)$ is an exact factorization of $\pgo$.

\subsubsection*{\it Restricted Lie algebras with triangular decompositions}
A triangular decomposition of a Lie algebra $\ago$ is a
collection $\left(\ago_{0}, \ago_{+}, \ago_{-},(\, |\, )\right)$ where 
$\ago_{0}, \ago_{-}, \ago_{+}$ are subalgebras of $\ago$ and $(\, |\, ): \ago \times \ago \rightarrow \ku$
is a  non-degenerate symmetric $\ago$-invariant bilinear form such that
$\ago_0$ is abelian, 
\begin{align*}
\ago &=\ago_{-} \oplus \ago_{0} \oplus \ago_{+}, &
\left[\ago_{\pm}, \ago_{0}\right] &\subset \ago_{\pm},& &\text{and} &
\left(\ago_{+} | \ago_{+}\right) =\left(\ago_{-} | \ago_{-}\right) &=\left(\ago_{+} | \ago_{0}\right)=\left(\ago_{0} | \ago_{-}\right) = 0. 
\end{align*}
A triangular decomposition gives rise to a Manin triple $(\pgo, \pgo_{1}, \pgo_{2})$
defined by
\begin{align*}
 \pgo &= \ago \times \ago,& \pgo_{1} &=\diag  \ago, & &\text{and}& \pgo_{2} &=\left\{\left(a_{-}+a_{0}, a_{+}-a_{0}\right): a_{\star} \in \ago_{\star}, \star \in \{+, 0, -\} \right\}. 
\end{align*}
If $\ago$ is restricted and $\ago_{0}, \ago_{-}, \ago_{+}$ are restricted subalgebras, then 
$\pgo$ is restricted and $(\pgo_{1}, \pgo_{2})$ is an exact factorization. There are other factorizations:
\begin{itemize}[leftmargin=*]
\item $(\ago_{\pm}, \ago_{0} \oplus \ago_{\mp})$ are exact factorizations of $\ago$;
\item $(\ago_{-} \oplus \ago_{0}, \ago_{0} \oplus \ago_{+})$ is an exact factorization of $\ago \times \ago_{0}$.
\end{itemize}

\subsubsection*{\it Restricted $\Z$-graded Lie algebras.}
A finite-dimensional $\Z$-graded Lie algebra $\s = \oplus_{i=-r}^{t} \s_{i}$ (where $r,t \in \N_0$)
is \emph{restricted} 
if the underlying Lie algebra $\s$ is restricted and 
$\s_{i}^{[p]} \subseteq \s_{pi}$ for all $i$. Then $\s_{+} \coloneqq  \oplus_{i=1}^{t} \s_{i}$
and $\s_{\leq 0} \coloneqq \oplus_{i=-r}^{0} \s_{i}$ form an exact factorization of $\s$. 

\section{Quasi-split extensions}\label{sect:quasi-split}
\subsection{Morita equivalence}
Two finite-dimensional  Hopf algebras 
$H$ and $U$ are \emph{Morita equivalent} ($H\mor U$) iff 
there exists an equivalence of braided tensor categories between the Drinfeld centers
$\mathcal Z(\rep H)$ and $\mathcal Z(\rep U)$, iff
 $D(H)$ and $D(U)$ are twist equivalent  quasitriangular Hopf 
 algebras.\footnote{Observe that this differs from 
  \cite{aapw}, where it was claimed  that $H\mor U \Leftrightarrow
 D(H) \simeq D(U)$  as quasitriangular Hopf algebras; we point out that this discrepancy  does 
 not affect the results of \emph{loc. cit.}}

In other words, $H\mor U$ if the tensor categories $\rep H$ and $\rep U$ are Morita equivalent. Notice that our defining condition is in fact a characterization of the original definition of Morita equivalence of tensor categories in \cite{muger}, 
\cite{EGNO}.

For instance, if either $U \simeq H^*$, or $U \simeq H^J$ is obtained from $H$
by twisting the comultiplication by $J \in H \otimes H$, 
or $U \simeq H_{\sigma}$ is obtained from $H$
by twisting the multiplication by $\sigma: H \otimes H \to \ku$,
then $H\mor U$.

\medbreak 
From now on, we assume that all Hopf algebras in \eqref{eq:exact} are finite-dimensional.
In this Section we study the following notion.

\begin{definition}\label{def:quasi-split}
We shall say that a cleft abelian exact sequence \eqref{eq:exact}  is \emph{quasi-split} if the Hopf algebra $C$ is Morita equivalent  to the bismash product $A \# B$ associated to the Singer pair arising from any cleaving maps.
\end{definition}

\subsection{Coquasi-Hopf algebras}
Recall that quasi-Hopf algebras were introduced by Drinfeld as generalizations of Hopf algebras, where the main difference is that the coassociativity of the comultiplication
holds up to a 3-tensor called the associator. Dually, a coquasi-bialgebra or 
coquasi-Hopf algebra $H$ is a generalization of a bialgebra or a Hopf algebra 
where the main difference is that the associativity of the multiplication
holds up to a dual 3-tensor $\varphi: H \otimes H \otimes H \to \ku$, 
called the coassociator. 
See e.g. \cite{schauenburg} for details.

\smallbreak
Here is the starting point of our analysis. 

\begin{definition}\cite{schauenburg}
Let $K$ and $Q$ be Hopf algebras. 
A \emph{generalized product coquasi-Hopf algebra} of $K$ and $Q$ 
is a co-quasi Hopf algebra $H$
together with coquasi-Hopf algebra maps 
\begin{align*}
i: K &\to H & &\text{and}& j: Q &\to H& &\text{ such that}& 
\operatorname{mult}(i \otimes j): K &\otimes Q \to H
\end{align*}
 is a linear isomorphism.
\end{definition}

\medbreak
Let us say that a finite-dimensional  Hopf algebra $L$
and a coquasi-Hopf algebra $U$ are \emph{Morita equivalent} iff 
there exists an equivalence of braided tensor categories between
$\mathcal Z(\corep L)$ and $\mathcal Z(\corep U)$ or equivalently that
the quantum doubles $D(L)$ and $D(U)$ are twist-equivalent.

\medbreak
By the results of \cite[Section 6]{schauenburg}, a cleft exact sequence of Hopf algebras \eqref{eq:exact} 
gives rise to a generalized product coquasi-Hopf algebra $H$ of $A^*$ and $B$, where 

\begin{itemize}[leftmargin=*]\renewcommand{\labelitemi}{$\circ$}
\item $H = B \otimes A^*$ as a coalgebra;

\smallbreak
\item  the multiplication is given by
\begin{align*}
(x \otimes \ell)(y \otimes m) &= x\left(\ell_{(1)} \fiz y_{(1)}\right) \otimes\left(\ell_{(2)} \trid y_{(2)}\right) m, & \ell,m &\in A^*, \ x, y \in B,
\end{align*}
where the maps $\fiz : A^* \otimes B \to A^*$ and $\trid : A^* \otimes B \to B$ are determined by the associated weak action $\rightharpoonup$ and the weak coaction $\rho$ by 
\begin{align*}
(\ell \fiz x) (a) &= \ell(x \rightharpoonup a),& 
\ell \trid x &= \rho(x)_i \, \ell(\rho(x)^i), & \ell &\in A^*,  x \in B;
\end{align*}

\smallbreak
\item the coassociator $\varphi: H \otimes H \otimes H \to \ku$ is determined by the cocycle $\sigma$ and the dual cocycle $\tau$ in the form
\begin{equation*}
\varphi(x \otimes \ell \otimes y \otimes m \otimes z \otimes r ) = \epsilon(x) \, \ell (y \rightharpoonup \tau^{(1)}(z_{(1)}) \sigma(y_{(2)}, \rho(z_{(2)})_i) \; m (\tau^{(1)}(z_{(1)}) \rho(z_{(2)})^i)) \, \epsilon (r).
\end{equation*} 
\end{itemize}

\begin{prop}\label{prop:cleft-morita} 
Given a cleft exact sequence of Hopf algebras \eqref{eq:exact},
 the Hopf algebra $C$ is Morita equivalent to the coquasi-Hopf algebra $H$. 
\end{prop}

\pf
The main result of \cite{schauenburg} implies the existence of an equivalence of monoidal categories 
\begin{equation*} {}_A(\corep C)_A \simeq  \corep H,
\end{equation*}
which amounts to the Morita equivalence of the categories $\corep C$, $\corep H$, 
and \emph{a fortiori} of $C$ and $H$. 
Indeed, an equivalence of braided tensor categories between the Drinfeld centers of $\corep C$ and $\corep H$ was established in \cite{schauenburg-jpaa}.
\epf

Combining Propositions \ref{prop:matched-pairs-hopf-classic}
and  \ref{prop:cleft-morita}, we obtain:

\begin{cor} Let $S$ be a finite-dimensional Hopf algebra and suppose that $S$ is a double crossproduct of its Hopf subalgebras $G$ and $L$. Then   $S$ is Morita equivalent to a bismash product $L^* \# G$. In particular, if $D(S)$ has fgc so does  $L^* \# G$. \qed
\end{cor}

\medbreak 
\subsection{Abelian extensions}
Suppose now that the exact sequence \eqref{eq:exact} is abelian. 
Then the coquasi-Hopf algebra $H$ described above turns out to be 
the double crossproduct associated to the Singer pair $(A, B)$ 
with a possibly nontrivial coassociator determined by the cocycles 
$\sigma$ and $\tau$. 
However, the coalgebra $H$ is cocommutative in this case. 

\begin{cor}\label{cor:quasi-split-phitrivial}
An abelian extension of $B$ by $A$ is quasi-split provided that the coassociator $\varphi$ is trivial. 
\end{cor}

\pf It follows from the fact that the split extension $A\# B$ is Morita equivalent to $A^* \bowtie B$ (with trivial associator) by Proposition \ref{prop:cleft-morita}. \epf

Exactness of the Kac sequence of \cite{schauenburg} implies that, since $A$ is 
finite-dimensional, $H$ is isomorphic as a  coquasi-Hopf algebra to the double crossproduct 
$A^*\bowtie B$ (as generalized products of $A^*$ and $B$) if and only if \eqref{eq:exact} is 
isomorphic (as a $B$-extension of $A$) to a twisting $(A \# B)_{ \chi}^J$ of the bismash 
product, where 
\begin{itemize}[leftmargin=*]\renewcommand{\labelitemi}{$\circ$}
\item $J \in A \otimes A$ is a twist in $A$, and

\item $\chi: B \otimes B$ is a 2-cocycle on  $B$, 
\end{itemize}
regarded respectively as a twist in $C$ and a 2-cocycle on $C$, see  
\cite[Theorem 6.3.6]{schauenburg}.

\medbreak
The next theorem is the main result of this section.

\begin{theorem} \label{thm:singer-pair}
Let $(A, B)$ be a Singer pair of finite-dimensional Hopf algebras. Given a quasi-split abelian extension $C$ of $B$ by $A$, the double $D(C)$ and a fortiori $C$ have fgc. 
In particular, $A\#B$ has fgc.
\end{theorem}

\pf We have that $C$ is Morita equivalent to the cocommutative Hopf algebra $A^* \bowtie B$. Whence  $D(C)$ is twist equivalent to $D(A^* \bowtie B)$, which has fgc by the main result of \cite{negron}. Hence $D(C)$ and therefore also $C$ have fgc.  \epf

\begin{example}\label{exa:mp-finite-gps}
Let $\varLambda$ and $\varGamma$ be  finite groups.
Consider an exact sequence of Hopf algebras
$\ku \to \ku ^{\varLambda} \to C \to \ku \varGamma \to \ku$ where $\ku ^{\varLambda}$ is the algebra of functions on $\varLambda$. Then $C \simeq \ku^{\varLambda}\#_{\sigma}^{\tau}\ku \varGamma$ is a bicrossed product. The relevant (weak) actions in this case are determined by actions by permutations $\trid:\varGamma \times \varLambda \to \varLambda$ and $\fiz: \varGamma \times \varLambda \to \varGamma$ that make  $(\varGamma, \varLambda, \trid, \fiz)$ into a matched pair of finite groups. Let $\varSigma = \varGamma \bowtie \varLambda$ be the associated double crossproduct group. 

\medbreak 
The Hopf algebra $C$ is Morita equivalent to a quasi-Hopf algebra $(\ku {\varSigma}, \omega)$, where $\omega \in H^3(\varSigma, \ku^\times)$ is the 3-cocycle attached to the class of $C$ under the Kac exact sequence (hence in particular, the restriction of $\omega$ to $\varGamma$ and $\varLambda$ is trivial). It follows from \cite{negron} that $C$ has fgc whenever $\omega$ is trivial. 

\medbreak 
For instance, we have $\ku \varSigma \mor \ku^{\varLambda}\#\ku \varGamma$, hence 
$\ku^{\varLambda}\#\ku \varGamma$ has fgc.
This is evident if $\car \ku$ is $0$ o coprime to $\vert \varSigma\vert$;
otherwise it follows alternatively
from \cite[Lemma 3.2.5]{aapw} since $\ku^{\varLambda}$ is semisimple.
\end{example}

\begin{example}\label{exa:mp-restricted-lie}
Let  $(\g, \lgo, \trid, \fiz)$ be a matched pair of restricted Lie algebras
and $\s = \g \bowtie \lgo$.
The corresponding matched pair of Hopf algebras 
$(\ugo(\g), \ugo(\lgo), \trid, \fiz)$ gives rise to an exact sequence 
$\ku \to \ugo(\lgo)^* \to \ugo(\lgo)^*\# \ugo(\g) \to \ugo(\g) \to \ku$.
We have $\ugo(\s)\mor \ugo(\lgo)^* \#\ugo(\g)$, hence 
$\ugo(\lgo)^* \#\ugo(\g)$ has fgc by \cite{negron}.
\end{example}

\section{The restricted Jordan plane}\label{sec:Jordan-plane}

\emph{In this Section and the next two, $\car \ku = p$ is an odd prime} (except when explicitly
stated otherwise).
In Section \ref{sec:hopf} we consider a subclass of the finite-dimensional
Nichols algebras introduced in \cite{aah-oddchar}
and show that their bosonizations with suitable group algebras are split abelian extensions,
therefore they have fgc.
In this Section and in the next, we work out the two simplest examples 
for illustration.

\subsection{The property fgc for bosonizations}
In this Subsection, $\car \ku$ is arbitrary.
We  record a useful result, a variation of \cite[Theorem 3.1.6]{aapw}.

\begin{theorem} \label{th:RtoRsmashH} 
Let $F$ be a finite group.

\begin{enumerate}[leftmargin=*,label=\rm{(\roman*)}]
\item\label{item:fgc-bos} Assume that $\ku F$ is semisimple. 
If $R$ is a finite-dimensional Hopf algebra in $\yd{\ku F}$ that has fgc, 
then $R\# \ku F$ has fgc.

\item\label{item:fgc-bos-dos} 
If $R$ is a finite-dimensional Hopf algebra in $\yd{\ku^F}$ that has fgc,
then $R\# \ku^F$ has fgc.
\end{enumerate} 
\end{theorem}

\pf  We sketch the proof for the reader's convenience. Let $K$ be either $\ku F$ as in \ref{item:fgc-bos} or $\ku^F$ as in \ref{item:fgc-bos-dos},
so clearly $K$ is semisimple. 
Let $R$ be as in  \ref{item:fgc-bos} or  \ref{item:fgc-bos-dos} accordingly.

Since the proofs of \cite[Corollary 3.13]{MPSW} and  \cite[Lemma 3.1.4]{aapw}
just require that $\ku$ is a field,  we conclude that the algebra
$\coh(R, \ku)$ is Noetherian. 
Now  \cite[Lemma 3.1.1]{aapw} also holds for any field, hence
$\coh(R, \ku)^{K}$ is finitely generated.

On the other hand, there is an isomorphism $\coh(R \# K, \ku) \simeq \coh(R, \ku)^{K}$, see 
\cite[Theorem 2.17]{stefan-vay}. Hence $\coh(R \# K, \ku)$ is finitely generated.

Next, given  a finitely generated $R\# K$-module $M$, one can prove that
$\coh(R\# K , M)$ is finitely generated as an
$\coh(R\# K , \ku)$-module repeating word-by-word the proof of the analogous fact in
\cite[Theorem 3.1.6]{aapw}. 
\epf

For further developments, it would be useful to extend Theorem \ref{th:RtoRsmashH} to
an arbitrary semisimple Hopf algebra $K$. 

\begin{lemma}\label{lema:nichols-fgc}
Let $A$ be a finite-dimensional Hopf algebra and $U \in \yd{A}$
such that $\toba(U)$ is finite-dimensional. If $\toba(U) \# A$ has fgc,
then so does the Nichols algebra $\toba(U)$. 
\end{lemma}

\pf
Since $\toba(U) \# A$ is free over $\toba(U)$,  \cite[Theorem 3.2.1]{ABFF}
implies the claim. 
\epf

\subsection{The minimal bosonization}\label{subsec:minimalbosonization}

We begin with the following basic example.
The \emph{block}  $\Vc(1,2)$ is the braided vector space
with a basis $\{x,y\}$ such that 
\begin{align}\label{equation:basis-block}
\begin{aligned}
c(x \ot  x) &=  x \ot  x,&c(y \ot  x) &=  x \ot  y, \\
c(x \ot  y) &=(y+ x) \ot  x,& c(y \ot  y) &=(y+x) \ot  y.
\end{aligned}
\end{align}

The Nichols algebra $\toba(\Vc(1,2))$ is
called the \emph{restricted Jordan plane}.

\begin{lemma} \label{lemma:restricted-jordan} \cite{clw}
The  restricted Jordan plane is generated by  $x$, $y$ with relations
\begin{align}\label{eq:rels B(V(1,2))}
&yx-xy+\tfrac{1}{2}x^2, &
&x^p, &
&y^p.
\end{align}
The set $\{ x^a y^b: 0 \le a,b< p\}$ is a basis of $\toba(\Vc(1,2))$, so $\dim \toba(\Vc(1,2)) = p^2$. \qed
\end{lemma}

The minimal bosonization of $\Vc(1,2)$ arises as follows.
Let $\Gamma = \Z/p = \langle g \rangle$.
We realize  $\Vc(1,2)$ in $\yd{\ku \Gamma}$ by 
\begin{align*}
g\cdot x &= x, & g\cdot y &= y + x,& \deg x = \deg y &= g.
\end{align*}

Thus  the Hopf algebra $H = \toba(\Vc(1,2)) \# \ku \Gamma$ has dimension $p^3$. 
We get a presentation of $H$  by generators $x, y, g$, where we 
identify $x= x\otimes 1, y = y \otimes 1, g = 1 \otimes g \in H$, 
with defining relations
\begin{align*}
x^p = y^p &= 0, & g^p &= 1,&  gx &= x g, & gy &= yg + xg, & 
yx &= xy - \tfrac{1}{2}x^2.
\end{align*}
The comultiplication of $H$ is determined by 
\begin{align*}
\Delta(g) &= g \otimes g, & \Delta(x)  &= x \otimes 1 + g \otimes x , 
& \Delta(y)  &= y \otimes 1 + g \otimes y .
\end{align*}
In addition, the monomials $g^ix^jy^\ell$, $0 \leq i, j, \ell \leq p-1$, form a basis of $H$.

\medbreak 
Let $K = \ku\langle x, g \rangle \subset H$ and $L = \ku[\zeta]/(\zeta^p)$
with $\zeta$ primitive. 

\begin{lemma} \label{lemma:restricted-jordan-ext} 
The Hopf algebra $H$  fits 
into a split abelian extension  $\ku \to K \overset{\iota}{\rightarrow} H \overset{\pi}{\rightarrow} L \to \ku$,
where $\iota$ is the inclusion and $\pi$ is defined by $\pi(x) = 0$, $\pi(g) = 1$
and $\pi(y) = \zeta$. 
\end{lemma}

\pf The defining relations of $H$ imply that $K$ is commutative. 
Furthermore $L$ is generated by the primitive element $\zeta$, whence cocommutative. 
Clearly $\pi$ is well-defined and $\ker \pi = HK^+$. 
We thus obtain an abelian exact sequence 
$\ku \to K \overset{\iota}{\rightarrow} H \overset{\pi}{\rightarrow} L \to \ku$.  

\smallbreak
Since $y^p = 0 = -\mathcal S(y)^p$, there exists a unique algebra map  $\mathpzc{s}: L \to H$ such that $\mathpzc{s}(\zeta) = -\mathcal S(y) = g^{-1}y$. 
Clearly, $\pi \mathpzc{s} = \id_L$. 
Being an algebra map, $\mathpzc{s}$ is automatically invertible for the convolution product. 
The $L$-colinearity of $\mathpzc{s}$ follows  from the relation 
\begin{align*}
(\mathpzc{s} \otimes \id) \Delta (\zeta) = g^{-1}y \otimes 1 + 1 \otimes \zeta = (\id \otimes \pi) \Delta \mathpzc{s}(\zeta) = (\id \otimes \pi) (g^{-1}y \otimes g^{-1} + 1 \otimes g^{-1}y).
\end{align*} 

Dually, let $\mathpzc{r}: H \to K$ be the linear map defined by 
\begin{align*}
\mathpzc{r}(g^i{x}^jy^\ell) &= g^i{x}^j  \,\, \text{ if } \,\, \ell  = 0, &
&\text{ and } &  \mathpzc{r}(g^i{x}^jy^\ell) &= 0 \text{ if } \,\, \ell > 0.
\end{align*} 
 It is clear that $\mathpzc{r}$ is $K$-linear, $\mathpzc{r}\iota = \id_K$ and $\mathpzc{r}\mathpzc{s} = \epsilon_L1_K$. 
 We next  show that $\mathpzc{r}$ is a coalgebra map. Since $\mathpzc{r}\vert_K = \id_K$, 
 we have $\Delta \mathpzc{r} (x)= (\mathpzc{r} \otimes \mathpzc{r} )\Delta(x)$, 
 for all $x\in K$.  
 Let $I = \ker \mathpzc{r} \subseteq H$, in other words, $I$ is the linear span of all monomials 
 $g^i{x}^jy^\ell$ with $\ell > 0$.
 The relation $\Delta(y) = y \otimes 1 + g \otimes y$ implies that 
 $\Delta(I) \subseteq I \otimes H + H \otimes I$. 
 Therefore $(\mathpzc{r} \otimes \mathpzc{r})\Delta(I) = 0$, 
 implying that $\mathpzc{r}$ is a coalgebra map. 
This shows that $(\mathpzc{s}, \mathpzc{r})$ is a splitting 
and finishes the proof of the lemma.
\epf

\begin{remark} \cite{A-Penha}
The Drinfeld double of $H$ fits into an abelian exact sequence
\begin{align*}
	\ku \to \mathbf R \rightarrow D(H)  
\rightarrow  {\mathfrak u}(\mathfrak{sl}_2(\ku)) \to \ku, 
\end{align*}
where  $\mathbf R$ is a local commutative Hopf algebra.
\end{remark}

The following result appeared already in \cite{NWW} with a different proof.

\begin{prop} \label{prop:restricted-jordan} 
The Hopf algebra $\toba(\Vc(1,2)) \# \ku \Gamma$ 
and the Nichols algebra $\toba(\Vc(1,2))$ have fgc.
\end{prop}
\pf
By  Lemma \ref{lemma:restricted-jordan-ext}, we may apply Theorem \ref{thm:fingencoh-extensions} and Lemma \ref{lema:nichols-fgc}.
\epf

\subsection{More bosonizations}

To deal with different realizations of the Jordan plane, we recall 
the notions of YD-pairs and YD-triples that are available in any characteristic. 
Let $A$ be a Hopf algebra. 

\medbreak
\begin{itemize}[leftmargin=*] \renewcommand{\labelitemi}{$\circ$}
\item A pair $(g, \chi) \in G(A) \times \Alg(A, \ku)$ is called a \emph{YD-pair} for $A$  if
\end{itemize}
\begin{align}\label{eq:yd-pair}
\chi(h)\,g  &= \chi(h\_{2}) h\_{1}\, g\, \Ss(h\_{3}),& h&\in A.
\end{align}

\medbreak
A YD-pair  $(g, \chi)$ gives rise to $\ku_g^{\chi} \in \yd{A}$ of dimension 1, with action 
and coaction  given by $\chi$ and $g$, respectively. Any one-dimensional object in
$\yd{A}$ is like this.

\medbreak
\begin{itemize}[leftmargin=*] \renewcommand{\labelitemi}{$\circ$}
\item A collection $(g, \chi, \eta)$ where
$(g, \chi)$ is a YD-pair for $A$ and $\eta \in \Der_{\chi,\chi}(A, \ku)$
is called a \emph{YD-triple} for $A$  if
\begin{align}\label{eq:YD-triple}
\eta(h) g &= \eta(h\_2) h\_1 g \Ss(h\_3), & h&\in A,
\\ \label{eq:YD-triple-jordan}
\text{and } \chi(g) &= \eta(g) = 1.
\end{align}

\end{itemize}
Notice that the existence of a YD-triple for $A$  when $\dim A < \infty$
forces that $\car \ku >0$, since $\eta \in A^*$ is a non-zero $(\chi, \chi)$- primitive.

\medbreak
A YD-triple $(g, \chi, \eta)$ gives rise to $\Vc_g(\chi,\eta) \in \yd{A}$, defined as the
vector space with a basis $\{x,y\}$, whose $A$-action and $A$-coaction are given by
\begin{align*}
h\cdot x &= \chi(h) x,& h\cdot y&=\chi(h) y + \eta(h)x,&h&\in A;& 
\delta(x) &= g\otimes x,& \delta(y) &= g\otimes y.
\end{align*}
By assumption \eqref{eq:YD-triple-jordan}, 
$\Vc_g(\chi, \eta)\simeq \Vc(1,2)$ as a braided vector space.

\medbreak
We now turn back to the assumption that $\car \ku = p > 2$.

\begin{question} Let $(g, \chi, \eta)$ be a YD-triple for $A$ where $\dim A < \infty$.
Does $\toba(\Vc_g(\chi, \eta)) \# A$ have fgc?
\end{question}

We point out that Theorem \ref{th:RtoRsmashH} 
does not apply in the present situation, by the following observation.

\begin{remark}
Let $A$ be a finite-dimensional Hopf algebra that admits a YD-triple $(g, \chi, \eta)$.
Then $A$ is not semisimple. 
\end{remark}

\pf Observe that the restriction $\eta:  \langle g\rangle \to \ku$ is a morphism 
of abelian groups, hence $\eta(g^p) = 0$, which implies that $p$ divides the order of $g$.
Thus $\ku \langle g\rangle$ is not semisimple, and a fortiori $A$ is not semisimple
by \cite[10.3.4]{radford-book}.
\epf

There are examples fitting into abelian exact sequences that are not necessarily split. 

\begin{remark}\label{rem:fgc-jordan-more-boson}
Let $F$ be a finite group  and  let $(g, \chi, \eta)$ be a YD-triple for $\ku F$.
It consists of $g \in Z(F)$, 
$\chi \in \widehat{F} \coloneqq \Hom_{\rm gps} (F, \ku^{\times})$ and
$\eta\in \Der_{\chi,\chi}(\ku F, \ku)$ such that $\chi(g) = \eta(g) = 1$.
Let $N \coloneqq \ker \chi \cap Z(F)\lhd F$. On one hand we consider the subalgebra
of $H  = \toba(\Vc_g(\chi, \eta)) \# \ku F$:
\begin{align*}
K &\coloneqq  \ku \langle x, \gamma: \gamma \in N\rangle 
\simeq \ku\langle x \rangle\# \ku N, 
\end{align*}

Here $K$ is commutative but not necessarily cocommutative.
On the other hand, $\chi$ induces a character $\overline{\chi}$ of $F/N$.
Let $\ku \zeta \in \yd{\ku (F/N)}$ corresponding to the YD-pair $(e, \overline{\chi})$.
Then 
\begin{align*}
L &\coloneqq \toba(\zeta) \# \ku (F/N) \simeq \ku[\zeta]/ (\zeta^p) \# \ku (F/N);
\end{align*}
here $\zeta$ is primitive. Clearly $L$ is cocommutative but not necessarily commutative.

\medbreak
Let $\iota: K \to H$ be the inclusion and let $\pi: H \to L$ be the map defined by $\pi(x) = 0$, $\pi(y) = \zeta$
and $\pi(\gamma) =$ the class of $\gamma$ in $F/N$. Then $H$ 
fits into the 
abelian exact sequence $\ku \to K \overset{\iota}{\rightarrow} H \overset{\pi}{\rightarrow} L \to \ku$, which is not split, for instance, when
the exact sequence of groups $1 \to N \rightarrow F\rightarrow F/N \to 1$
 is not split.
\end{remark}

In the setting of the previous remark, the Hopf algebra $H^*$ has fgc.

\begin{prop}
Let $H  = \toba(\Vc_g(\chi, \eta)) \# \ku F$, where $F$ is a finite group and $(g, \chi, \eta)$ is a YD-triple for $\ku F$. Then $H^*$ has fgc.
\end{prop}

\pf Arguing as in \cite[Lemma 1.5]{A-Penha}, we see that $H^* \simeq \Rc \# \ku^{F}$,
where $\Rc \simeq \toba(W)$ is isomorphic to $\toba(\Vs(1,2))$ as algebras (although not as braided Hopf algebras). Anyway, $\Rc$ has fgc by Proposition \ref{prop:restricted-jordan}, hence $H^*$ has fgc by  Theorem \ref{th:RtoRsmashH}.
\epf

\section{The first Laestrygonian algebra \texorpdfstring{$\toba(\lstr_q( 1,\ghost))$}{}} \label{section:lstr-11disc}
\subsection{The Nichols algebra \texorpdfstring{$\toba(\lstr_q( 1,\ghost))$}{}} \label{subsection:lstr-11disc}
The next  example of our interest depends on the 
data: $q \in \kut$ and $a \in {\Fp}^{\times}$.  Let 
\begin{align*}
\sa &\in \{1- p, 2-p, \dots,-2, -1 \} & &\text{such that }&
\sa &\equiv 2a \mod p. 
\end{align*}
The \emph{ghost} is the  integer $\ghost\coloneqq  -\sa \in \I_{p-1}$; 
since $p$ is odd, $\ghost$ determines $a$.
To this data we attach the braided vector space
 $\lstr_q( 1, \ghost)$  with  basis $\mathpzc{b} = \{x_1, y_1, x_2\}$ and  braiding given by
\begin{align}\label{eq:braiding-block-point}
(c(b\otimes b'))_{b, b'\in \mathpzc{b}} &=
\begin{pmatrix}
 x_1 \otimes x_1&  ( y_1 + x_1) \otimes x_1& q \,  x_2  \otimes x_1
\\
 x_1 \otimes y_1 & ( y_1 + x_1) \otimes y_1&  q \, x_2  \otimes y_1
\\
q^{-1} x_1 \otimes x_2 &  q^{-1}(y_1 + a x_1) \otimes x_2&   x_2  \otimes x_2
\end{pmatrix}.
\end{align}
Thus $V_1 \coloneqq \ku x_1 + \ku y_1  \simeq \Vc(1, 2)$ and $V_2 \coloneqq \ku x_2 $ 
satisfy $c: V_i \otimes V_j = V_j \otimes V_i$, $i,j \in \{1,2\}$; in particular 
$V_1$ and $V_2$ are braided subspaces of $V$.
We introduce 
\begin{align}\label{eq:def-zn}
z_0 &\coloneqq x_2, &
z_{n+1} & \coloneqq y_1z_n - qz_ny_1,& n &>0.
\end{align}

\begin{lemma} \label{lemma:lstr-11disc} \cite[4.3.1]{aah-oddchar}
The algebra
$\toba(\lstr_q( 1, \ghost))$ is presented by generators $x_1,y_1, x_2$ and relations \eqref{eq:rels B(V(1,2))}, and, in the notation \eqref{eq:def-zn},
\begin{align}
x_1x_2&=q \, x_2x_1,  \label{eq:lstr-rels&11disc-1} \\
z_{1+\ghost}&=0,  \label{eq:lstr-rels&11disc-qserre} \\
z_tz_{t+1}&=q^{-1} \, z_{t+1}z_t, &  0\le & t < \ghost, \label{eq:lstr-rels&11disc-2}
\\
z_t^p&=0, & 0\le & t \le  \ghost. \label{eq:lstr-rels&11disc-pot-p}\end{align}
 The algebra $\toba(\lstr_q( 1, \ghost))$ has a PBW-basis
\begin{equation}\label{base}B=\{ x_1^{m_1} y_1^{m_2} z_{\ghost}^{n_{\ghost}} \dots z_1^{n_1} z_0^{n_0}: m_i, n_j \in \I_{0, p}\},
\end{equation} hence $\dim \toba(\lstr_q( 1, \ghost)) = p^{\ghost + 3}$.
\qed \end{lemma}

\subsection{ A suitable realization} \label{subsection:lstr-11disc-realization}
In order to realize $\lstr_q( 1, \ghost)$ in $\yd{\ku \Gamma}$ for some finite group $\Gamma$, we need to assume that $q$ is a root of 1. Set $d \coloneqq \ord q$; 
then $d$ is coprime to $p = \car \ku$. Fix a positive integer $f$ which is a multiple of $pd$.
A suitable choice of group is 
\begin{align*}
\Gamma = {\Z/f} \times {\Z/f}  &= \langle g_1 \rangle \oplus \langle g_2 \rangle,&
\text{where } \ord g_1 &= \ord g_2 = f.
\end{align*}
It is not difficult to see that  $\lstr_q( 1, \ghost)$ can be realized in $\yd{\ku \Gamma}$ by 
\begin{align}\label{eq:realization-lestr}
\begin{aligned}
g_1\cdot x_1 &= x_1,& g_1\cdot y_1 &= y_1 + x_1,& g_1\cdot x_2 &= q x_2,
\\
g_2\cdot x_1 &= q^{-1}x_1,& g_2\cdot y_1 &= q^{-1}(y_1 + a x_1),& g_2\cdot x_2 &= x_2,
\\
\deg x_1 &= g_1, &\deg y_1 &= g_1, &\deg x_2 &= g_2.
\end{aligned}
\end{align}

\medbreak
Therefore  the Hopf algebra $\Huno \coloneqq \toba(\lstr_q( 1, \ghost))\# \ku \Gamma$
 is presented by generators $x_1$, $y_1$, $x_2$, $g_1$, $g_2$, with relations \eqref{eq:rels B(V(1,2))}, \eqref{eq:lstr-rels&11disc-1}, \eqref{eq:lstr-rels&11disc-qserre}, \eqref{eq:lstr-rels&11disc-2}, \eqref{eq:lstr-rels&11disc-pot-p}, and
\begin{align}
g_1^{f} &= 1,&   g_2^{f} &= 1,&  g_1g_2 &= g_2g_1,  
\\
g_1x_1 &= x_1g_1,&  g_1y_1 & = y_1g_1 + g_1x_1, &  g_1x_2 &= qx_2g_1,  
 \\
g_2x_1 &= q^{-1}x_1g_2, & g_2y_1 &= q^{-1}(y_1g_2 + ax_1g_2),&  g_2x_2 &= x_2g_2.
 \end{align}
 The comultiplication of $\Huno$ is determined by $\Delta(g_i) = g_i \otimes g_i$, $i = 1, 2$, and
 \begin{equation}\label{deltaxi}\Delta(x_1) = x_1 \otimes 1 + g_1 \otimes x_1, \quad \Delta(y_1) = y_1 \otimes 1 + g_1 \otimes y_1, \quad 
 \Delta(x_2) = x_2 \otimes 1 + g_2 \otimes x_2.
\end{equation}
Clearly $\dim \Huno = p^{\ghost + 3} f^2$.

\subsection{The split case} \label{subsection:lstr-11disc-q=1}
In this subsection we deal with $\toba(\lstr_1( 1, \ghost))$; Let $\Gamma$ be as above
with $f$ divisible by $pd$ with $d = \ord q$.
Let $K \subseteq H \coloneqq \toba(\lstr_1( 1, \ghost))\# \ku \Gamma$ be the Hopf subalgebra  
\begin{align*}
K = \ku\langle x_1, g_1, g_2 \rangle. 
\end{align*}
We shall consider the restricted enveloping algebra 
\begin{align*}
L = \ugo(\lgo),
\end{align*}
where $\lgo$ is the restricted Lie algebra defined as follows.
Let  $V(\ghost)$ be the simple $\mathfrak{sl}(2)$-module of highest weight $\ghost$.
Let $E = (\begin{smallmatrix}0 &1 \\ 0& 0 \end{smallmatrix})$. Pick
a basis $(v_n)_{n\in \I_{0, \ghost}}$ of $V(\ghost)$ such that
\begin{align*}
E \cdot v_n &= v_{n+1}, & n&\in \I_{0, \ghost - 1}, &E \cdot v_\ghost &= 0.
\end{align*}
Let $\lgo =  V(\ghost)\rtimes \ku E$, a Lie subalgebra of the motion Lie algebra  
$V(\ghost)\rtimes \mathfrak{sl}(2)$;  it follows from Lemma \ref{lemma:mp-restricted} 
that $\lgo$ is  restricted with $p$-operation equal to $0$. 
That is, $L$ is presented by generators $v_0$ and $E$
with defining relations, in terms of $ v_{n+1} = Ev_n - v_nE$, $ n\in \I_{0, \ghost - 1}$,
\begin{align*}
v_mv_n &= v_n v_m, &v_n^p &= 0, &  m, n&\in \I_{0, \ghost}, &Ev_\ghost - v_\ghost E&= 0,
& E^p &= 0.
\end{align*}

\begin{prop} \label{prop:restricted-lstr-ext} 
The Hopf algebra $H= \toba(\lstr_1( 1, \ghost))\# \ku \Gamma$ fits into 
a split abelian exact sequence  $\ku \to  K\overset{\iota}\to H \overset{\pi}\to L \to \ku$,
where $\iota$ is the inclusion and $\pi$ is defined by 
\begin{align}\label{eq:restricted-lstr-ext}
\pi(g_1) &= 1,&  \pi(g_2) &= 1, &  \pi(x_1) &= 0, &  \pi(y_1) & = E, &  \pi(x_2)& = v_0.
\end{align}
\end{prop}

\pf By assumption, the Hopf subalgebra $K$ is commutative. In addition $L$ is cocommutative. Notice that $\dim K = pf^{2}$. 
We see by inspection that $K$ is normal, i.e., $HK^+ = K^+H$. Thus we have an abelian exact sequence 
\begin{align*}
\ku \to K\overset{\iota}\to H \to H/ HK^+\to \ku.
\end{align*}
Since $\dim H = p^{\ghost + 3} f^2$, $\dim H/ HK^+ = p^{\ghost + 2} = \dim L$.
Now \eqref{eq:restricted-lstr-ext} determines an algebra map $\pi: H \to L$ 
by \eqref{eq:def-zn} and the defining relations of $\toba(\lstr_1( 1, \ghost))$. 
The map $\pi$ is surjective,  has
$HK^+ \subseteq \ker \pi$ and preserves the comultiplication 
since the classes of $y_1$ and $x_2$ are primitive in $H/ HK^+$.
By dimension counting, $\pi$ induces an isomorphism of Hopf algebras $H/ HK^+ \to  L$. 

\medbreak 
Observe that $\pi(z_n) = v_n$, for all $n \geq 0$.
To define a splitting, we argue as in the proof of  Lemma \ref{lemma:restricted-jordan-ext}. The universal property of $\ugo(\lgo)$ implies the existence of a unique algebra map 
$\mathpzc{s}: L = \ugo(\lgo) \to H$ such that 
$$\mathpzc{s}(E) = -\mathcal S(y_1), \quad \mathpzc{s}(v_n) = -\mathcal S(z_n), \quad \forall n = 0, 1, \dots$$
Clearly $\pi \mathpzc{s} = \id_L$. 
To see that $\mathpzc{s}$ is $L$-colinear, it will be enough to verify that the maps 
$(\mathpzc{s} \otimes \id_L)\Delta$ and $(\id \otimes \pi)\Delta \mathpzc{s}$ 
agree on $E$ and $v_0$, 
because both $(\mathpzc{s} \otimes \id_L)\Delta$  and $(\id_H \otimes \pi)\Delta \mathpzc{s}$ 
are algebra maps and $E$ and $v_0$ generate $L$ as an algebra. 
This follows at once from \eqref{deltaxi}.

\medbreak
Let now $\mathpzc{r}: H \to K$ be the linear map defined on the basis of $H$ 
arising from \eqref{base} by 
$$\mathpzc{r}(g_1^cg_2^bx_1^{m_1} y_1^{m_2} z_{\ghost}^{n_{\ghost}} \dots z_1^{n_1} z_0^{n_0}) 
= \begin{cases} g_1^cg_2^bx_1^{m_1}, \; \textrm{if } m_2 + n_{\ghost} + \dots + n_0 = 0,\\  0, \quad \textrm{otherwise}\end{cases}$$ 
As in the proof of Lemma  \ref{lemma:restricted-jordan-ext}, 
we see that $\mathpzc{r}$ is a $K$-linear retraction of $\iota$ and that $\mathpzc{r}\mathpzc{s} = \epsilon_L 1_K$. 
Moreover, $\mathpzc{r}$ is a coalgebra map. 
To see this, we consider the subspace $I = \ker \mathpzc{r}$ which coincides 
with the linear span of all monomials   
$g_1^cg_2^bx_1^{m_1} y_1^{m_2} z_{\ghost}^{n_{\ghost}} \dots z_1^{n_1} z_0^{n_0}$ 
such that $m_2 + n_{\ghost} + \dots + n_0 > 0$. 
The defining relations of $H$ imply that $I$ is a left ideal. 
We claim that $\Delta(I) \subseteq  I \otimes H + H \otimes I$. 
Indeed, \eqref{deltaxi} implies that $\Delta(y_1), \Delta(x_2) \in I \otimes H + H \otimes I$. 
 By \cite[4.2.10, pp. 28]{aah-triang} and the comultiplication formula of the bosonization, 
 there exist $\nu_{j, n} \in \ku$ such that
\begin{align*} 
\Delta(z_n) &= z_n \otimes 1 + \sum_{j = 0}^n \nu_{j, n} \; x_1^{n-j} g_1^j g_2 \otimes z_j,
& \forall n &\in \I_{0, \ghost}.
\end{align*}

Thus $\Delta(z_n) \in I \otimes H + H \otimes I$, for all $n \geq 0$. Then  the claim follows. 
Since $\mathpzc{r}\vert_K = \id_K$, we  conclude that $\mathpzc{r}$ is a coalgebra map.
Thus $(\mathpzc{s}, \mathpzc{r})$ is a splitting and the proof of the Proposition is complete.
\epf

\begin{prop} \label{prop:restricted-lstr-ext-fgc} 
The Hopf algebra $H = \toba(\lstr_1( 1, \ghost))\# \ku \Gamma$ 
and its double $D(H)$ have fgc.
\end{prop}
\pf
This follows from Proposition \ref{prop:restricted-lstr-ext} and Theorem \ref{thm:singer-pair}.
\epf

\subsection{The general case} \label{subsection:lstr-11disc-general}
Recall that $q \in \ku^\times$ has  order $d$, that $f$ is a multiple of $pd$ and that
  $\Gamma = \langle g_1, g_2: \, g_1^f = g_2^f = 1, \, g_1g_2 = g_2g_1 \rangle 
\simeq \Z/f \times \Z/f$.

\begin{theorem} \label{thm:restricted-lstr-ext} 
The Hopf algebra $\Huno = \toba(\lstr_q( 1, \ghost))\# \ku \Gamma$ has fgc.
\end{theorem}
 \pf Consider the bilinear form $\sigma: \Gamma \times \Gamma \to \ku^\times$ determined by   
\begin{equation*} \sigma(g_i, g_j) = 
\begin{cases}
1, \quad i = j,\\
q, \quad i = 1, j = 2,\\
1, \quad i = 2, j = 1.
\end{cases}
\end{equation*}
Then $\sigma$ is a 2-cocycle in $\Gamma$.  Let $\vartheta: \Gamma \times \Gamma \to \ku^\times$, $\vartheta(g, h) = \sigma(g, h)\sigma(h, g)^{-1}$ be the associated antisymmetric bilinear form. Thus, 
\begin{equation*} \vartheta(g_i, g_j) = 
\begin{cases}
1, \quad i = j,\\
q, \quad i = 1, j = 2,\\
q^{-1}, \quad i = 2, j = 1.
\end{cases}
\end{equation*}

 As in Subsection \ref{subsubsec:twist} of the Appendix, let $\Fc_\sigma: \yd{\ku\Gamma} \to \yd{\ku\Gamma}$ be the  monoidal functor associated to $\sigma$. 
The image under $\Fc_\sigma$ of the  braided vector space
 $V = \lstr_1( 1, \ghost)$ in \ref{subsection:lstr-11disc-realization} is  $\Fc_{\sigma}(V) = V$
 with the  same grading as $V$ and $\Gamma$-action \eqref{eq:twistingaction}; that is, 
 \begin{align}\label{eq:realization-lestr-twisted}
\begin{aligned}
g_1\cdot x_1 &= x_1,& g_1\cdot y_1 &= y_1 + x_1,& g_1\cdot x_2 &= q x_2,
\\
g_2\cdot x_1 &= q^{-1}x_1,& g_2\cdot y_1 &= q^{-1}(y_1 + a x_1),& g_2\cdot x_2 &= x_2,
\\
\deg x_1 &= g_1, &\deg y_1 &= g_1, &\deg x_2 &= g_2.
\end{aligned}
\end{align}
Hence $\Fc_\sigma(V) = \lstr_q( 1, \ghost)$. By Lemma \ref{lem:twisting-Nichols}. the bosonization  $\Huno = \toba(\lstr_q( 1, \ghost))\# \ku \Gamma$ is a cocycle deformation of the bosonization $H = \toba(\lstr_1( 1, \ghost))\# \ku \Gamma$. 
By Proposition \ref{prop:restricted-lstr-ext-fgc}, the double 
$D\left(H\right)$ has fgc.
Hence $\Huno$  also has fgc, as claimed.
 \epf

\begin{cor} \label{cor:restricted-lstr-ext} Let $q \in \G_{\infty}$.
The Nichols algebra $\toba(\lstr_q( 1, \ghost))$ has fgc. \qed
\end{cor}

\section{Pointed Hopf algebras over abelian groups}\label{sec:hopf}

\subsection{A class of braided vector spaces}
We now proceed with a family of Hopf algebras which are bosonizations of some  Nichols 
algebras introduced in \cite{aah-oddchar}, analogues in odd characteristic  of  
Nichols algebras appearing  in \cite{aah-triang}.
We  show that, up to a cocycle deformation, 
the Hopf algebra $H = \toba(\Vs(\bq, \ba))\# \ku \Gamma$
fits into a split abelian exact sequence 
$\ku \to K \rightarrow H \rightarrow \ugo(\lgo) \to \ku$,
where the restricted Lie algebra $\lgo$ is determined explicitly.

\medbreak
\subsubsection*{Data}
We shall consider braided vector spaces depending on:

\medbreak
\begin{itemize}  [leftmargin=*]\renewcommand{\labelitemi}{$\circ$}
\item two positive integers $t < \theta$;

\medbreak
\item a matrix $\bq = (q_{ij})_{i,j \in \I_{\theta}}$ such that 
\begin{align}\label{eq:matrix-q}
q_{ij}q_{ji} &= 1, & q_{ii} &= 1, & i, j &\in \I_{\theta},\,  i \neq  j.
\end{align}

\medbreak
\item a family  $\ba = (a_{ij}) _{\substack{i \in \I_{t + 1, \theta}, \\  j \in  \I_{t}}}$
with entries in $\fp$. We lift this family to $\Z$ as follows:  
\end{itemize}

\begin{itemize}\renewcommand{\labelitemi}{$\diamond$}
\item when $a_{ij} \neq 0$, we take
\begin{align}\label{eq:discrete-ghost}
\sa_{ij} &\in \{1- p, 2-p, \dots,-2, -1 \} & &\text{such that }&
\sa_{ij} &\equiv 2a_{ij} \mod p,
\end{align}
and then we set $\ghost_{i,j}\coloneqq  -\sa_{ij} \in \I_{p-1}$.

\medbreak
\item when $a_{ij} = 0$, we lift it to $\ghost_{i,j} = 0$.
\end{itemize}

\medbreak
The family $\ghost \coloneqq (\ghost_{i,j}) _{\substack{i \in \I_{t + 1, \theta}, \\  j \in  \I_{t}}}$, equivalent to   $\ba$, is  the \emph{ghost};  both
$\ba$ and $\ghost$ are needed.

\medbreak
\subsubsection*{Definition}
We define a braided vector space $(\Vs(\bq, \ba), c)$ from the data above.
First, it has a decomposition 
$\Vs(\bq, \ba) =  V_{1} \oplus \dots \oplus V_t \oplus \dots \oplus V_\theta$
such that
\begin{align*}
c(V_i \otimes V_j) &=  V_j \otimes V_i, \quad i,j \in \I_{\theta}.
\end{align*}

Then we assume:
\begin{itemize}
\item If $j\in \I_t$, then $V_j \simeq \Vc(1,2)$ (the blocks). 
Let  $\{x_j, y_j\}$ be a basis of $V_j$ realizing
\eqref{equation:basis-block}:
\begin{align}\label{equation:basis-block-i}
\begin{aligned}
c(x_j \ot  x_j) &=  x_j \ot  x_j,&c(y_j \ot  x_j) &=  x_j \ot  y_j, \\
c(x_j \ot  y_j) &=(y_j +  x_j) \ot  x_j, & c(y_j \ot  y_j) &=(y_j + x_j) \ot  y_j.
\end{aligned}
\end{align}

\bigbreak
\item If $i\in \I_{t + 1, \theta}$, then $\dim V_{i} = 1$; these are the points. We fix a basis $\{x_i\}$ of $V_i$. Thus, $\{x_i: i \in \I_{\theta}\} \amalg \{y_j: j \in \I_{t}\}$ is a basis of $\Vs(\bq, \ba)$.
\end{itemize}

\medbreak 
The braidings $c_{ij} = c_{\vert V_i\otimes V_j}$, $i,j\in \I_{\theta}$ are given
by the data $\bq$ and $\ba$ as follows:

\begin{itemize}[leftmargin=*] \renewcommand{\labelitemi}{$\circ$}
\item If $i ,j\in \I_{\theta}$, then 
\begin{align}\label{eq:braiding-qij}
c(x_i \otimes x_j) &= q_{ij}  x_j \otimes x_i.
\end{align}

\medbreak
\item If $i ,j\in \I_{t}$ are blocks, then $c_{ij} = q_{ij} \tau$, $\tau$ being the  flip,
for $i \neq j$ while $c_{jj}$ is given by \eqref{equation:basis-block-i}.

\medbreak
\item If $j\in \I_t$ a block and $i \in \I_{t+1, \theta}$  is a point, 
then  $c_{\vert (V_i \oplus V_j)\otimes (V_i \oplus V_j)}$ is given  by \eqref{eq:braiding-qij}
and 
\begin{align}\label{eq:intro-braiding-block-point}
c( y_j \otimes x_i) &= q_{ji} x_i  \otimes y_j, &
c( x_i \otimes y_j) &= q_{ij} (y_j + a_{ij} x_j )  \otimes x_i.
\end{align}
\end{itemize}

\medbreak
As in \cite{aah-triang},  $\Vs(\bq, \ba)$ is described by a diagram of  the following shape:

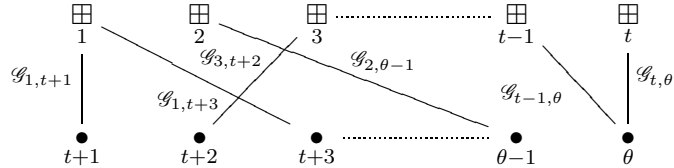
\begin{align*}
\xymatrix{\underset{1}{\boxplus}   \ar @{-}[d]_{\ghost_{1, t+1}}
\ar @{-}[rrd]_<(0.5){ \ghost_{1, t+3} }  & \underset{2}{\boxplus} \ar@{-}[rrrd]^<(0.5){\ghost_{2, \theta -1}} 
& \underset{3}{\boxplus} \ar  @{.}[rr] \ar @{-}[ld]_<(0.4){\ghost_{3,t+2}\hspace{-5pt}} & & \underset{t-1}{\boxplus}  
\ar @{-}[rd]_{\ghost_{t-1, \theta}}
& \underset{t}{\boxplus} \ar @{-}[d]^{\ghost_{t, \theta}}   
\\ \underset{t+1}{\bullet}   & \underset{t+2}{\bullet}   
& \underset{t+3}{\bullet} 
\ar  @{.}[rr]
&&\underset{\theta -1}{\bullet}   
& \underset{\theta}{\bullet}}
\end{align*}
That is, there are $\mathpzc{r}$ blocks, $\theta - t$ points and a line decorated by $\ghost_{k, \ell}$
when  $\ghost_{k, \ell} \neq 0$, joining the block $k$ with the point $\ell$;
 this graph is admissible in the sense of \cite[Definition 1.3.7]{aah-triang}.

\begin{remark} The data above is an $\ab$-triple $\triple = (\bn, \bq, \ba)$
of rank $\theta$, cf. \S \, \ref{subsec:appendix-bvs-ab}, as follows:
\begin{itemize} 
\item $\bn = (n_j)_{j\in \I_{\theta}}$ is  given by $n_j = 2$, if $1 \le j \le t$;
and $n_j = 1$, if $t <  j \le \theta$.

\medbreak
\item the matrix $\bq = (q_{ij})_{i,j \in \I_{\theta}}$ is as given above
with the constraint \eqref{eq:matrix-q}; 

\medbreak \item the family $\tb = (\tb_{ij}) _{i,j\in \I_{\theta,}}$  
where $\tb_{ij}\in \End \ku^{n_j}$ is given by
\begin{align*}
\tb_{ij} &= 0, &&&\text{when } t  &<  j \le \theta \text{ or } 1 \le i, j \le t; 
\\
\tb_{ij}(x_{j}) &= 0, &\tb_{ij}(y_j) &= a_{ij} x_{j},
&\text{when }  1 &\le  j \le t \text{ and }  t <  i \le \theta.
\end{align*}

\end{itemize}

If $t= \theta = 1$, then we have the restricted Jordan plane as in Section \ref{sec:Jordan-plane}. 
If $t=1$ and $\theta = 2$, then we recover $\lstr_q( 1, \ghost)$ as in Section \ref{section:lstr-11disc}.
\end{remark}

\subsubsection*{Grading} The subspace $\mathscr{U}$ 
of $\Vs = \Vs(\bq, \ba)$ spanned by $(x_i)_{i \in \I_{\theta}}$ is a braided subspace
and $ \mathscr{U} \subseteq \Vs$ is a filtration of braided vector spaces.
The associated graded vector space $\gr \Vs = \mathscr{U} \oplus \Vs/ \mathscr{U}$
is of diagonal type; indeed, if $\overline{y}_j$ is the class of $y_j$ 
in $\Vs/ \mathscr{U}$, $j\in \I_{t} $, then
$(x_i)_{i \in \I_{\theta}}  \amalg (\overline{y}_j)_{ j \in \I_{t}}$ 
is a basis of $\gr \Vs$ and the braiding of $\gr \Vs$ is given by
\begin{align*}
&\begin{aligned}
c(x_i \ot  x_k) &=  q_{ik} \, x_k \ot  x_j,&c(\overline{y}_j \ot  x_i) &= q_{ji}  \, x_i \ot  \overline{y}_j, \\
c(x_i \ot  \overline{y}_j) &= q_{ij} \, \overline{y}_j \ot  x_i,
& c(\overline{y}_j \ot  \overline{y}_{\ell}) &= q_{j \ell} \,  \overline{y}_{\ell} \ot  \overline{y}_j,
\end{aligned}&
\text{for all } i,k &\in \I_{\theta}, \  j, \ell \in \I_{t}.
\end{align*}
Precisely, let  $Q \simeq \Z^{\theta + t}$ be the free abelian group with basis
$(\alpha_i)_{i \in \I_{\theta}} \amalg (\beta_j)_{ j \in \I_{t}}$; then the braiding
is given by the bilinear form $\bp: Q \times Q \to \ku^{\times}$   defined by
\begin{align*}
\bp(\alpha_i \ot  \alpha_j)  = \bp(\alpha_i \ot  \beta_j) 
= \bp(\beta_i \ot  \alpha_j) = \bp(\beta_i \ot  \beta_j)   &=  q_{ij}
\end{align*}
whenever $i$ and $j$ make sense. By \eqref{eq:matrix-q}, $\toba(\gr \Vs)$ 
is a quantum linear space and we have
\begin{align*}
\bp(\gamma \ot  \delta)  &= \bp(\delta \ot  \gamma)& \forall \gamma, \delta    &\in Q.
\end{align*}
Now consider the grading of $\Vs$ and its extension to $T(\Vs)$
given by $\deg x_i = \alpha_{i}$, 
$\deg y_j = \beta_{j}$. 
\medbreak
Given  points $h, \ell \in \I_{t+1, \theta} $ we set
\begin{align}
\label{eq:defn-sch-several-blocks}
\sch_{h, \bn} &\coloneqq  (\ad_c y_1)^{n_1} \dots (\ad_c y_t)^{n_{t}} x_{h}, 
\qquad \bn = (n_1,\dots,n_{t})\in\N^t_0;
\\
\label{eq:Ah}
\mathcal A_h &\coloneqq  \{\bn \in\N^t_0: 0 \le \bn \le \ghost_h = (\ghost_{h,1}, \dots, \ghost_{h,t})\},
\text{ ordered lexicographically.}
\end{align}
In terms of the bilinear form $\bp$, we also set
\begin{align*}
\bp_{h, \ell;\bm,\bn} &\coloneqq  \bp_{\,\deg \sch_{h, \bm}, \deg \sch_{\ell, \bn}}
= \big(\prod_{k,j\in\I_t} q_{kj}^{m_kn_j}\big)\big(\prod_{k\in\I_t} q_{k \ell}^{m_k}\big) \big(\prod_{j\in\I_t} q_{h  j}^{n_j}\big) q_{h \ell}, 
\quad \bm \in \Ac_h, \ \bn \in \Ac_{\ell}.
\end{align*}

\begin{remark}\label{remark:sch-adjoint}
In the definition \eqref{eq:defn-sch-several-blocks}, $\ad_c$ means the braided adjoint.
However, it could be replaced by a sequence of $q$-commutators (with various $q$). 
More precisely,
given $h \in \I_{t+1, \theta}$ and  $\bn = (n_1,\dots,n_{t})\in\N^t_0$,
set $\sch = \sch_{h, \bn}$. Then
\begin{align}\label{eq:sch-conmutation}
\sch &= y_j \, \widetilde{\sch} - q \, \widetilde{\sch} \, y_j,
\end{align}
where $j = \min \{i\in \I_t: n_i \neq 0\}$, 
$\widetilde{\bn} = (0,\dots, 0, n_j-1, n_{j+1}, \dots, n_{t})$,
$\widetilde{\sch} = \sch_{h, \widetilde{\bn}}$ and
\begin{align*}
q &= q_{jh} \prod_{j \le i \le t} q_{ji}^{n_i}.
\end{align*}
See the proof of \cite[Lemma 7.2.3]{aah-triang}. In particular,
if $\bq = \uno$ is the matrix with all entries equal to $1$, then
$q = 1$ and we can replace $\ad_c$ by the usual adjoint in \eqref{eq:defn-sch-several-blocks}.
\end{remark}

\subsubsection*{The Nichols algebra $\toba(\Vs(\bq, \ba))$}
The following result is not included in \cite{aah-oddchar} but its proof is
similar to that of \cite[4.3.1, 6.1]{aah-oddchar}.

\begin{lemma} \label{lemma:general}
The algebra $\toba(\Vs(\bq, \ba))$ is presented by generators 
$x_i$, $i\in \I_{\theta}$, $y_j$, $j\in \I_{t}  = \I_{t}$, and relations
\begin{align}\label{eq:poseidon-defrels-Jordan}
&x_j^p =0, \quad y_j^p =0, \quad
y_jx_j -x_jy_j+\tfrac{1}{2}x_j^2 =0, & j &\in\I_{t} ;
\\\label{eq:poseidon-defrels-blocks-commute}
&x_kx_j = q_{kj} \, x_jx_k, \, x_ky_j = q_{kj} \, y_jx_k, \quad y_ky_j = q_{kj} \, y_jy_k
& k \neq j &\in\I_{t} ;
\\ \label{eq:poseidon-defrels-q-commute}
&x_jx_{h} = q_{j h} \, x_{h} x_j, & j & \in\I_{t} , \,\, h \in \I_{t+1, \theta} ;
\\ \label{eq:poseidon-defrels-q-Serre}
&(\ad_c y_j)^{1+\ghost_{jh}}(x_{h})=0, & j& \in\I_{t} , \,\, h \in \I_{t+1, \theta} ;
\\ \label{eq:poseidon-rels-K-3}
&\sch_{h, \bn}^p =0, & \bn &\in \Ac_h, \,\, h \in \I_{t+1, \theta} 
\\
\label{eq:poseidon-rels-K-1}
&\sch_{h, \bm}\sch_{\ell, \bn} = \bp_{h, \ell;\bm,\bn} \, \sch_{\ell, \bn}\sch_{h, \bm}, 
& h, \ell \in \I_{t+1, \theta} , \,     &\bm \in \Ac_h, \,  \bn \in \Ac_\ell.
\end{align}
A basis of $\toba(\Vs(\bq, \ba))$ is given by
\begin{align}\label{eq:base-general}
B =\big\{ x_1^{ m_1} y_1^{m_2} \dots x_t^{ m_{2t-1}} y_t^{m_{2t}} 
\prod_{\substack{h\in \I_{t+1, \theta}\\ \bn \in \Ac_h}} \sch_{h, \bn}^{n_{h,\bn}}: \, 
& 0\le n_{h, \bn}, m_j  < p \mbox{ if }  j\in \I_{2t},
h\in \I_{t+1, \theta},  \bn \in \Ac_h \big\}.
\end{align}
Hence $\dim \toba(\Vs(\bq, \ba)) =  p^{2t + \sum_{h\in \I_{t+1, \theta}} \vert \Ac_h\vert}$.
\end{lemma}
\pf
Argue as in the proofs of \cite[4.3.1, 6.1]{aah-oddchar}; relation 
\eqref{eq:poseidon-rels-K-1} follows as in the discussion in 
\cite[pp. 107]{aah-triang}.
\epf 

\subsection{A suitable realization}
In order to realize $\Vs(\bq, \ba)$ in $\yd{\ku \Gamma}$ for some finite group $\Gamma$, 
we  need to assume that the entries of the matrix $\bq$ are roots of 1.  Set
\begin{align*}
d \coloneqq \lcm\{\ord q _{ij} :i, j \in \I_{\theta}\};
\end{align*}
then $d$ is coprime to $p = \car \ku$. Fix a positive integer $f$ which is a multiple of $pd$.
A suitable choice of group is  
\begin{align*}
\Gamma &= ({\Z/f}) ^{\theta} = \langle g_1 \rangle \oplus \cdots \oplus \langle g_{\theta}\rangle, &
&\text{where}& \ord g_1 &= \ord g_2 = \dots \ord g_{\theta} = f.
\end{align*}
In other words, $\Gamma$ is generated by $g_1, \dots, g_{\theta}$ with relations
\begin{align}\label{eq:rels-Gamma}
g_i^f  &=  1,& g_ig_j &= g_jg_i,&  i, j &\in \I_{\theta}.
\end{align}
It is not difficult to see that  $\Vs(\bq, \ba)$ can be realized in $\yd{\ku \Gamma}$ by 
\begin{align}\label{eq:realization-V(q,a)}
\begin{aligned}
g_k\cdot V_\ell &= q_{k \ell} \id_{V_{\ell}},& k&\in \I_{t}, \, \ell\in \I_{\theta}, \, k \neq \ell; 
\text{ or } k, \ell \in \I_{t+1, \theta};
\\
g_i\cdot x_j  &= q_{ij} x_j, \quad
g_i\cdot y_j = q_{ij} (y_j + a_{ij} x_j ),& i&\in \I_{t+1, \theta}, j\in \I_{t},
\\
g_j\cdot x_j  &=  x_j, \quad
g_j\cdot y_j =   y_j + x_j,&  j &\in \I_{t},
\\
\deg x_\ell &= g_{ \ell}, \quad \deg y_j = g_{j}, & \ell &\in \I_{\theta}.  \, j \in \I_{t}.
\end{aligned}
\end{align}

Let $\Huno$ denote the Hopf algebra  $\toba(\Vs(\bq, \ba))\# \ku \Gamma$;
it has a presentation by generators $x_i$, $y_j$, $g_k$, $i,k \in \I_{\theta}$,
$j\in \I_{t}$ with  relations
\eqref{eq:poseidon-defrels-Jordan}, \eqref{eq:poseidon-defrels-blocks-commute},
\eqref{eq:poseidon-defrels-q-commute}, \eqref{eq:poseidon-defrels-q-Serre}, 
\eqref{eq:poseidon-rels-K-3}, \eqref{eq:poseidon-rels-K-1}, 
\eqref{eq:rels-Gamma} and those induced by \eqref{eq:realization-V(q,a)}.
The comultiplication is given by
\begin{align}\label{eq:comult-H-gen}
\Delta(x_i) &= x_i \otimes 1 + g_i \otimes x_i, &
\Delta(y_j) &= y_j \otimes 1 + g_j \otimes y_j, &
\Delta(g_i) &= g_i \otimes g_i, \, 
 i\in \I_{\theta},  j\in \I_{t}.
\end{align}
One has
\begin{align*}
\dim \Huno =  p^{2t + \sum_{h\in \I_{t+1, \theta}} \vert \Ac_h\vert}f^{\theta}.
\end{align*}

\subsection{The split case} \label{subsection:V(1,a)}
In this subsection we deal with $\toba(\Vs(\uno, \ba))$, where $\uno \in \ku^{\theta}$
has all entries equal to $1$.
Let $\Gamma$ be as above with $f$ divisible by $pd$ with $d = \ord q$.
Consider  the   Hopf subalgebra  $K \subseteq H = \toba(\Vs(\uno, \ba))\# \ku \Gamma$ 
\begin{align*}
K = \ku\langle x_1, \dots, x_{t}, \, g_1, \dots, g_{\theta} \rangle. 
\end{align*}

We shall consider the restricted enveloping algebra 
\begin{align*}
L = \ugo(\lgo),
\end{align*}
where $\lgo$ is the restricted Lie algebra defined by the  following steps.

\begin{itemize}  [leftmargin=*]\renewcommand{\labelitemi}{$\circ$}

\item $ \g \coloneqq  \g_1\oplus \cdots \oplus \g_t$ is the direct sum of $\mathpzc{r}$ copies 
$\g_j \simeq\mathfrak{sl}(2)$, $j \in \I_{t}$.

\medbreak
\item $E_j \in \g_j$ is the element corresponding to 
$(\begin{smallmatrix}0 &1 \\ 0& 0 \end{smallmatrix})$; 
$\ngo \coloneqq \ku E_1\oplus \cdots \oplus \ku E_t \hookrightarrow \g$.

\medbreak 
\item
For any point $h \in \I_{t + 1, \theta}$, let $V(\ghost_{h,j})$ be the simple $\g_j$-module of highest weight $\ghost_{h,j}$.
Pick a basis $(v_{h,j; n})_{n\in \I_{0, \ghost_{h, j}}}$ of $V(\ghost_{h,j})$ such that
\begin{align*}
E_j \cdot v_{h,j; m} &= v_{h, j; m+1}, & 0 \leq m &< \ghost_{h,j}, 
&E_j \cdot v_{h,\ghost_{hij}} &= 0.
\end{align*}

\medbreak 
\item For any point $h \in \I_{t + 1, \theta}$, let $V(\ghost_{h})$ be the simple $\g$-module 
\begin{align*}
V(\ghost_{h}) \coloneqq V(\ghost_{h,1}) \otimes \cdots \otimes V(\ghost_{h,t}).
\end{align*}
Recalling the notation \eqref{eq:Ah}, a basis of $V(\ghost_{h})$ is formed by the elements
\begin{align*}
v_{h, \bm } &\coloneqq  v_{h,1, [m_1]} \otimes \cdots \otimes v_{h, t, [m_t]}, &
\bm  &=  (m_1,\dots,m_{t}) \in \Ac_h.
\end{align*}
For $h \in \I_{t + 1, \theta}$, set $v_{h} \coloneqq  v_{h,1, [0]} \otimes \cdots \otimes v_{h, t, [0]}$. 
Then for any $\bm \in \Ac_h$, we have
\begin{align*}
v_{h, \bm } &=  E_1^{m_1} \cdots E_t^{m_t} \cdot  v_{h}.
\end{align*}
Indeed,  $V(\ghost_{h})$ is the simple $\g$-module of highest weight 
$\ghost_{h} = (\ghost_{h, 1},\dots, \ghost_{h, t})$.

\medbreak
\item  Finally, let  $\lgo =  V(\ghost)\rtimes \ngo$, where 
\begin{align*}
V(\ghost) \coloneqq V(\ghost_{t+1}) \oplus \cdots \oplus V(\ghost_{\theta}).
\end{align*}
\end{itemize} 
By Lemma \ref{lemma:mp-restricted},  $\lgo$ is  a restricted Lie algebra 
with $p$-operation equal to $0$.  
The restricted enveloping algebra  $L = \ugo(\lgo)$ 
is presented by generators $E_j$, $j\in \I_t$, and $v_{h}$, $h\in \I_{t+1, \theta}$;  set 
\begin{align*}
v_{h, \bm } &=  (\ad E_1)^{m_1} \cdots (\ad{E_t})^{m_t}  (v_{h}).
\end{align*}
Then the defining relations are 
\begin{align}\label{eq:lie2}
E_jE_k &= E_kE_j, &   j, k  &\in\I_t;
\\ \label{eq:lie4}
(\ad E_{j})^{1+\ghost_{jh}}(v_{h}) &= 0, & j& \in\I_t, \,\, h \in \I_{t+1, \theta};
\\
\label{eq:lie5}
v_{h, \bm } v_{i, \bn } &= v_{i, \bn } v_{h, \bm }, & i,h &\in \I_{t+1, \theta}, \,\,
 \bm  \in \Ac_h,  \,\,\bn \in \Ac_i;
\\ 
\label{eq:lie1}
E_j^p &=0, &  j &\in\I_t;
\\
\label{eq:lie6}
v_{h, \bm }^p &=0, & h &\in \I_{t+1, \theta}, \,\, \bm \in \Ac_h.
\end{align}
One has
\begin{align*}
	\dim L =  p^{t + \sum_{h\in \I_{t+1, \theta}} \vert \Ac_h\vert}.
\end{align*}
\begin{remark} \label{rem:lie} 
Clearly $\lgo$ is a Lie subalgebra of $V(\ghost)\rtimes \g$; 
we do not need the reference to $\g$ here but it will be necessary in further developments.
\end{remark}

\begin{prop} \label{prop:restricted-gen-ext} 
The Hopf algebra $H = \toba(\Vs(\uno, \ba))\# \ku \Gamma$  fits 
into a split abelian exact sequence  
$\ku \to K \overset{\iota}{\rightarrow} H \overset{\pi}{\rightarrow} L \to \ku$,
where $\iota$ is the inclusion and $\pi: H \to L$ is determined by 
\begin{align}\label{eq:restricted-gen-ext}
\pi(g_i) &= 1,& i \in \I_\theta,&  &   \pi(x_j) &= 0, &  \pi(y_j) & = E_j, & j \in \I_t, &  & \pi(x_\ell)& = v_\ell, &\ell \in \I_{t+1, \theta}&.
\end{align}
\end{prop}

The proof has the same pattern as the proof of Proposition \ref{prop:restricted-lstr-ext}.

\pf Evidently $K$ is commutative, $L$ is cocommutative and $\dim K = p^tf^{\theta}$. 
By inspection, $K$ is normal, i.e., $HK^+ = K^+H$. Thus we have an abelian exact sequence 
\begin{align*}
	\ku \to K\overset{\iota}\to H \to H/ HK^+\to \ku.
\end{align*}
We have $\dim H/ HK^+ = \dim L$.
By the defining relations of $\toba(\Vs(\uno, \ba))$, see Remark \ref{remark:sch-adjoint},
the assignment \eqref{eq:restricted-gen-ext} determines a Hopf algebra map $\pi: H \to L$,
which  is surjective and has
$HK^+ \subseteq \ker \pi$. Thus $\pi$ induces an isomorphism of Hopf algebras $H/ HK^+ \to  L$. 

\medbreak 
Observe that $\pi(\sch_{h, \bn}) = v_{h, \bn}$, for all $h, \bn$.
The section $\mathpzc{s}: L = \ugo(\lgo) \to H$ is the unique algebra map  such that 
\begin{align*}
\mathpzc{s}(E_j) &= -\mathcal S(y_j) = -g_j^{-1} y_j,
 & \mathpzc{s}(v_h) &= -\mathcal S(x_h) = -g_h^{-1} x_h, &
 \forall  j &\in \I_{t}, \, h\in \I_{t + 1, \theta}.
\end{align*}
We appeal again to Remark \ref{remark:sch-adjoint} to see that $\mathpzc{s}$ is indeed well-defined. 
As  in the proof of Proposition \ref{prop:restricted-lstr-ext}, we
get that $\mathpzc{s}$ is an $L$-colinear section of $\pi$.

\medbreak
Observe that the generic element of the basis  of $H$ arising from \eqref{eq:base-general}
can be expressed as 
\begin{align*}
u = g_1^{b_1}\dots g_\theta^{b_{\theta}} x_1^{ m_1}  \dots x_t^{ m_{t}} 
y_1^{n_1} \dots y_t^{n_{t}} 
\prod_{\substack{h\in \I_{t+1, \theta}\\ \bn \in \Ac_h}} \sch_{h, \bn}^{n_{h, \bn}}.
\end{align*}
Let  $\mathpzc{r}: H \to K$ be the linear map defined on $u$ by the formula
$$\mathpzc{r}(u)
= \begin{cases} g_1^{b_1}\dots g_\theta^{b_{\theta}} x_1^{ m_1}  \dots x_t^{ m_{t}}, \; \textrm{if } n_{k}  =  n_{h, \bn} = 0, \,  \forall k,h,\bn;
\\  0, \quad \textrm{otherwise.}\end{cases}$$ 

Clearly, $\mathpzc{r}$ is a $K$-linear retraction of $\iota$ and $\mathpzc{r}\mathpzc{s} = \epsilon_L 1_K$. 
We claim that $\mathpzc{r}$ is a coalgebra map. 
To see this, we consider the subspace $I = \ker \mathpzc{r}$ which is 
 the linear span of the monomials   
 \begin{align*}
 g_1^{b_1}\dots g_\theta^{b_{\theta}} x_1^{ m_1}  \dots x_t^{ m_{t}} 
 y_1^{n_1} \dots y_t^{n_{t}} 
 \prod_{\substack{h\in \I_{t+1, \theta}\\ \bn \in \Ac_h}} \sch_{h, \bn}^{n_{h, \bn}}: \qquad
 n_1 + \dots + n_h + \sum_{\substack{h\in \I_{t+1, \theta}\\ \bn \in \Ac_h}} n_{h, \bn} > 0.
 \end{align*}
By the defining relations of $H$, $I$ is a left ideal. 

We claim that $\Delta(I) \subseteq  I \otimes H + H \otimes I$. First,  
$\Delta(y_j), \Delta(x_h) \in I \otimes H + H \otimes I$
for all $j \in \I_{t}$, $h\in \I_{t + 1, \theta}$ by \eqref{eq:comult-H-gen}. 
Next, by \cite[Lemmas 7.2.3 (b) and 7.2.4]{aah-triang} 
and the comultiplication formula of the bosonization, 
we have
\begin{align*} 
\Delta(\sch_{h, \bn} ) &\in \sch_{h, \bn}  \otimes 1 + 
\sum_{0 \le \mathbf{k} \le \bn } H \otimes \sch_{h, \mathbf{k}} 
\subset I \otimes H + H \otimes I,
& \forall h &, \bn.
\end{align*}
Let $u$, $v$ be monomials in $I$ such that
$\Delta(u), \Delta(v) \in I \otimes H + H \otimes I$. 
Then $\Delta(u v) \in I \otimes H + H \otimes I$. By a recursive argument,
the second claim follows. 
Since $\mathpzc{r}\vert_K = \id_K$, we  get the first claim, i.e.,
 that $\mathpzc{r}$ is a coalgebra map.
Thus $(\mathpzc{s}, \mathpzc{r})$ is a splitting and the proof of the Proposition is complete.
\epf

\begin{prop} \label{prop:restricted-lstr-ext-V(q,a)} 
The Hopf algebra $H = \toba(\Vs(\uno, \ba))\# \ku \Gamma$ 
and its double $D(H)$ have fgc.
\end{prop}
\pf
This follows from Proposition \ref{prop:restricted-gen-ext} and Theorem \ref{thm:singer-pair}.
\epf

\subsection{The general case} \label{subsection:V(q,a)-general}

\begin{theorem} \label{thm:restricted-gen-ext} 
The Hopf algebra $\Huno = \toba(\Vs(\bq, \ba))\# \ku \Gamma$ has fgc.
\end{theorem}

\pf
By Corollary \ref{cor:twist-equivalent}, 
the bosonization  $\Huno = \toba(\Vs(\bq, \ba))\# \ku \Gamma$ 
is a cocycle deformation of the bosonization $H = \toba(\Vs(\uno, \ba))\# \ku \Gamma$. 
By Proposition \ref{prop:restricted-lstr-ext-V(q,a)}, the double 
$D\left(H\right)$ has fgc. Hence $\Huno$  also has fgc, as claimed. 
\epf

Since $H$ is free over $\toba(\Vs(\bq, \ba))$, 
 Theorem \ref{thm:restricted-gen-ext} together with  \cite[Theorem 3.2.1]{ABFF}
 implies:
\begin{cor}
The Nichols algebra $\toba(\Vs(\bq, \ba))$ has fgc. \qed
\end{cor}

 \appendix
 \section{Cocycle-equivalence of Nichols algebras}

 In this Appendix we describe braided vector spaces over abelian groups
 and spell out conditions for their Nichols algebras being twist-equivalent.
 Except in Corollary \ref{cor:twist-equivalent}, $\car \ku$ is arbitrary.
 
 \subsection{Braided vector spaces over abelian groups} \label{subsec:appendix-bvs-ab}
 
 \begin{definition} 
Let $\theta \in \N$.
An  $\ab$-triple (of rank $\theta$) is a collection $\triple = (\bn, \bq, \tb)$ where
\begin{itemize} 
\item $\bn = (n_j)_{j\in \I_{\theta}}$ is  a family of positive integers,  normalized
by $n_1 \geq n_2 \geq \dots \geq n_{\theta}$;

\medbreak
\item $\bq = (q_{ij})_{i,j \in \I_{\theta}}$ is a matrix with invertible entries and 

\medbreak \item $\tb = (\tb_{ij}) _{i,j\in \I_{\theta,}}$ is a family where $\tb_{ij}\in \End \ku^{n_j}$
 satisfies $\tb_{ij} = 0$ when $\dim V_j = 1$ and $\tb_{ik}\tb_{jk} = \tb_{jk}\tb_{ik}$ for all $i,j,k$.
 \end{itemize}
 An $\ab$-triple is \emph{nilpotent} if every $\tb_{ij}$ is nilpotent.
  \end{definition}
 
 We attach a braided vector space to an $\ab$-triple $(\bn, \bq, \tb)$
 by the following recipe.
 Let $V = \oplus_{j\in \I_{\theta}} V_{j}$ be a vector space  with a decomposition
such that $\dim V_j = n_j$ for all $j\in \I_{\theta}$;  
pick a basis of $V_j$, pull back $\tb_{ij}$ 
 to $\tb_{ij}\in \End V_j$ and define $c \in GL(V\otimes V)$  by 
 \begin{align}\label{eq:braiding-abelian}
 c(x\otimes y) &= q_{ij}\left(y + \tb_{ij}(y) \right) \otimes x, & x\in V_i, \quad  y\in V_j.
 \end{align}
 
 The proof of the following result is left to the reader.
 
 \begin{lemma} \label{lema:realization-Lambda} 
 	 The pair  $(V,c)$ is a braided vector space that can be realized 
 	 over $\yd{\ku \Lambda}$, where 
 $\Lambda\simeq \Z^\theta$ with canonical  basis 
 $\alpha_{1}, \dots, \alpha_{\theta}$, by
 \begin{align*}
 V_{\alpha_i} &= V_i, & \alpha_i \rightharpoonup x &= q_{ij}\left(x + \tb_{ij}(x) \right),
 & x\in V_j, \ i,j\in \I_{\theta}. \qed 
 \end{align*}
 \end{lemma}

 \subsection{Cocycle equivalence} 
We start by discussing the relationship between cocycle deformation and bosonization, as it appears in \cite{majid-oeckl}.
Let $H$ be  a Hopf algebra;  
let  $\sigma : H \otimes H \to \ku$ be an invertible 2-cocycle;
 let $H_{\sigma}$ be the Hopf algebra which is $H$ as coalgebra
and has  multiplication  
 \begin{align}\label{eq:twist-mult}
 x \cdot_{\sigma }y &= \sigma(x\_{1}, y\_{1})x\_{2}y\_{2}\sigma^{-1}(x\_{3}, y\_{3}), &x,y &\in H,
 \end{align}
 Let now $R$ be a Hopf algebra in $\yd{H}$ and $A \coloneqq R \#H$ the bosonization 
  with canonical projection and injection
 $\pi : A \to H$ and $\iota: H \to A$. Let $\sigma^\pi : A \otimes  A \to \ku$ be given by $\sigma^\pi \coloneqq\sigma(\pi \otimes \pi)$;
 this is an invertible 2-cocycle  on $A$. 
 Then $\pi : A_{\sigma^{\pi}} \to H_{\sigma}$ and $\iota: H_{\sigma} \to A_{\sigma^{\pi}}$ are still Hopf algebra maps. Hence $A_{\sigma^{\pi}} \simeq R_{\sigma} \# H_{\sigma}$
 where $R_{\sigma}$  is a Hopf algebra in $\yd{H_{\sigma}}$ that coincides with $R$ as vector subspace of $A$, with multiplication
 \begin{align}\label{eq:twist-mult-R}
 x\cdot_{\sigma } y &= \sigma(x\_{0}, y\_{0})x\_{1}y\_{1}, & x,y &\in R_{\sigma}.
 \end{align} 
 
 \begin{lemma} \label{lema:twisting-nichols}
 \cite[2.13]{A-Schneider-cambr} 
 If $R = \oplus_{n\in \N_0}R(n)$ is a  graded Hopf algebra in $\yd{H}$, then
 $R_{\sigma}$ is a  graded Hopf algebra in $\yd{H_{\sigma}}$
 with $R(n)=R_{\sigma}(n)$ as vector spaces for all $n \geq 0$. 
 Also, $R$ is a Nichols algebra if and only if $R_\sigma$ is. \qed
 \end{lemma}

 \subsection{Cocycles over an abelian group}\label{subsubsec:twist}
 We fix an abelian group $\Gamma$. 
 Let   $\sigma: \Gamma \times \Gamma \to \ku^{\times}$ be a group  2-cocycle:
 $\sigma(gh,k) \sigma(g,h) = \sigma(g,hk) \sigma(h,k)$, for all $g,h,k\in \Gamma$.
 Then the map $\vartheta:  \Gamma \times \Gamma \to \ku^{\times}$  given by
 \begin{align*}
 \vartheta(g,h) &= \sigma(g,h) \sigma^{-1}(h,g), & g,h &\in \Gamma,
 \end{align*}
 is  bilinear and antisymmetric.
 Let $\Fc_{\sigma}: \yd{\ku\Gamma} \to \yd{\ku\Gamma}$ be the monoidal functor that assigns $V \mapsto \Fc_{\sigma}(V) = V$
 with the same grading and action
 \begin{align}\label{eq:twistingaction}
 g\rightharpoonup_{\sigma} v &=   \vartheta(g, h) g \rightharpoonup v, & g, h\in \Gamma, v \in V_h,
 \end{align}
 where $\rightharpoonup$ is the action on $V$.
 The braiding $c_{\sigma}$ in  $\Fc_{\sigma}(V) \otimes \Fc_{\sigma}(W)$ is given  by
 \begin{align} \label{twistingbraiding}
 c_{\sigma} (v \otimes w) &= \vartheta(g,h) c (v \otimes w),&
 v\in V_{g}, w&\in W_h, \ g, h \in \Gamma.
 \end{align}
 
 Fix $V\in \yd{\Gamma}$ and let $A = \toba(V)\# \ku \Gamma$; as before $\pi: A \to \ku \Gamma$ and $\iota: \ku \Gamma \to A$ are the canonical projection and the inclusion. 
 Clearly, the linear extension $\sigma: \ku \Gamma \otimes \ku \Gamma \to \ku$ of $\sigma$
 is a 2-cocycle as in the previous subsection.
 Let $\sigma^\pi = \sigma (\pi \otimes \pi)$ be as above. 
 
 \begin{lemma}\label{lem:twisting-Nichols} The Hopf algebra $A_{\sigma^\pi}$ is  isomorphic to $\toba(V)_{\sigma}  \# \ku \Gamma$ 
 with  the same comultiplication of $\toba(V)$ and the multiplication  given by
 \begin{align}
 \label{eq:twistingalgebra}
 x._{\sigma^\pi}y &= \sigma(g,h) xy, & x&\in \toba(V)_g, \  y \in \toba(V)_h, & g,h &\in \Gamma.
 \end{align}
 Furthermore $\toba(V)_{\sigma}  \simeq \toba(\Fc_{\sigma}(V))$.
 \end{lemma}
 
 \noindent\emph{Proof.} The first part follows from the discussion above; \eqref{eq:twistingalgebra} is particular instance of \eqref{eq:twist-mult-R}.
 For the second part we show that the degree one homogenous component of $\toba(V)_{\sigma}$ is $\Fc_{\sigma}(V)$
 and apply Lemma \ref{lema:twisting-nichols}.
 If $y\in V_h$, then
 \begin{align*} 
 g._{\sigma^\pi} y &=  \sigma(g, \pi(y)) g  \sigma^{-1}(g, 1)
 + \sigma(g, h) \, gy\sigma^{-1}(g, 1) + \sigma(g, h)\, gh \sigma^{-1}(g, \pi(y))
 \\ &=   \sigma(g, h) gy;
 \end{align*}
 Similarly, $(gy)._{\sigma^\pi} g^{-1} = \sigma(gh, g^{-1} ) gyg^{-1}  \sigma^{-1}(g,g^{-1} )$ so the action is given by
 \begin{align*} 
 (g._{\sigma^\pi} y) ._{\sigma^\pi} g^{-1}  
 &=   \sigma(g, h) (gy._{\sigma^\pi} g^{-1})   =   \sigma(g, h)\sigma(gh, g^{-1} ) gyg^{-1}  \sigma^{-1}(g,g^{-1} )
 \\& =    \sigma(g, h)\sigma(h, g)^{-1} g\rightharpoonup y = g\rightharpoonup_{\sigma^\pi}y. \hspace{90pt}\qed
 \end{align*}
 
\begin{definition}  \cite[\S 2.4]{ABFF,A-Schneider-cambr}.
 Two braided Hopf algebras $R$ and $S$ are \emph{cocycle-equivalent}
 if there exist a Hopf algebra $H$ and an invertible 2-cocycle $\sigma : H \otimes H \to \ku$ such that
 \begin{itemize}[leftmargin=*]
 \item $R$ is realizable in $\yd{H}$; 
 
 \item $S$ is isomorphic to $R_{\sigma}$ as a braided Hopf algebra.
 \end{itemize}
\end{definition}

 The following definition extends \cite[Lemma 4.3]{ABFF}
 and \cite[\S 2.4]{A-Schneider-cambr}.
 
 \begin{definition}\label{def:twist-equivalent} 
Two braided vector spaces $(V, c)$ and $(V', c')$ arising from $\ab$-triples $(\bn, \bq, \tb)$ and $(\bn', \bq', \tb')$ 
 are {\it twist-equivalent} 
 if   
 \begin{align}\label{eq:twistequiv}
 \begin{aligned}
 \bn  &= \bn', & \tb &= \tb',
&
 q_{ii} &= q'_{ii}, & q_{ij}q_{ji} &= q'_{ij}q'_{ji}, &  i, j &\in \I_{\theta}.
 \end{aligned}\end{align} 
 \end{definition}

 \begin{lemma}\label{prop:twist-equivalent}
 If  the braided vector spaces $(V, c)$ and $(V', c')$ arising from the 
 $\ab$-triples $(\bn, \bq, \tb)$ and $(\bn', \bq', \tb')$  are twist-equivalent, then  
 the Nichols algebras $\toba(V)$ and $\toba(V')$ are cocycle-equivalent.
 \end{lemma}
 
 \pf
 By \eqref{eq:twistequiv} there exists a linear isomorphism $\psi: V \to V'$
 preserving the decompositions $V = \oplus_{i\in \I_{\theta}} V_i$ and $V' = \oplus_{i\in \I_{\theta}} V'_i$
 and intertwining the endomorphisms  $\tb_{ij}$ and $\tb'_{ij}$. We realize $(V, c)$ as in Lemma \ref{lema:realization-Lambda}.
 We  consider the unique bilinear form  $\sigma: \Lambda \times \Lambda \to \ku^{\times}$, hence a group 2-cocycle,
 given by
 \begin{equation}\label{eq:appendix-cocycle}
 \sigma(\alpha_{i}, \alpha_{j}) = \begin{cases} &q'_{ij} q_{ij}^{-1} , \qquad i\le j, \\
 &1, \qquad \qquad i > j.\end{cases}  \end{equation} 
 
 We claim that $\psi: V_{\sigma} \to V'$ 
 is an isomorphism in $\yd{\ku \Lambda}$.  Clearly $\psi$ preserves the grading.
 Assume $i\le j$ and let $x\in V_i$, $y\in V_j$. Then
 \begin{align*}
 \alpha_{j}\rightharpoonup_{\sigma}x & \overset{\eqref{eq:twistingaction}}{=}\sigma(\alpha_{j}, \alpha_{i})\sigma^{-1}(\alpha_{i}, \alpha_{j})q_{ji}\left(x + \tb_{ji}(x) \right) 
 \\ &= (q'_{ij})^{-1}  q_{ij} q_{ji} \left(x + \tb_{ji}(x) \right) 
 \overset{\eqref{eq:twistequiv}}{=} q'_{ji} \left(x + \tb_{ji}(x) \right);
 \\
 \alpha_{i}\rightharpoonup_{\sigma}y & \overset{\eqref{eq:twistingaction}}{=} 
 \sigma(\alpha_{i}, \alpha_{j})\sigma^{-1}(\alpha_{j}, \alpha_{i})q_{ij} \left(y + \tb_{ij}(i) \right) 
 = q'_{ij} \left(y + \tb_{ij}(i) \right) .
 \end{align*}
 Thus  $\psi$ preserves the action of $\Lambda$, hence it  extends to an isomorphism 
 $\Psi: \toba (V_{\sigma}) \to \toba(V')$ of Hopf algebras in $\yd{\ku \Lambda}$. 
 Now $\toba (\Fc_{\sigma}(V)) \simeq \toba(V)_{\sigma}$  by Lemma \ref{lem:twisting-Nichols}.
 \epf
 
The following statement is needed in the paper. Assume that $\car \ku = p$ is odd.
 Let $\Vs(\bq, \ba)$ and $\Vs(\bq', \ba)$ be two braided vector spaces
 as in Section \ref{sec:hopf} with the same $\theta$.
By hypothesis, $ q_{ii} = q'_{ii} = 1$, and $ q_{ij}q_{ji} = q'_{ij}q'_{ji} = 1$, 
for all $ i, j \in \I_{\theta}$.
 
\begin{cor}\label{cor:twist-equivalent} 
Assume that there exists a positive integer $f$  such that
\begin{align}\label{eq:cor-app-p}
p &\text{ divides } f,
\\ \label{eq:cor-app-qij}
\ord q_{ij} &\text{ divides } f, & \ord q'_{ij} &\text{ divides } f, & \forall i,j &\in \I_{\theta}.
\end{align}
Let  $\Gamma = (\Z/f)^{\theta}$. 
 Then  $\Vs(\bq, \ba)$ and $\Vs(\bq', \ba)$ are realizable in $\yd{\ku \Gamma}$
 and there exists  an invertible 2-cocycle 
 $\sigma :  \ku \Gamma \otimes \ku \Gamma \to \ku$ such that
 $\toba(\Vs(\bq', \ba))$ is isomorphic to $\toba(\Vs(\bq, \ba))_{\sigma}$ as Hopf algebras in $\yd{\ku \Gamma}$.
\end{cor}
 
 \pf The proof of the first claim on the realizations is straightforward using the hypotheses
 \eqref{eq:cor-app-p}  and \eqref{eq:cor-app-qij}.
 By \eqref{eq:cor-app-qij}, there is a
 unique bilinear form  $\sigma: \Gamma \times \Gamma \to \ku^{\times}$
 given by \eqref{eq:appendix-cocycle}.
 Then the second claim follows as in the proof of Lemma \ref{prop:twist-equivalent}.
  \epf


\begin{thebibliography}{NWW}

\bibitem{andersen-jantzen} H. H. Andersen and J. C.  Jantzen.
\emph{Cohomology of induced representations for algebraic groups}.
Math. Ann. \textbf{269}, pp. 487--525 (1984).

\bibitem{andrus-leyva} N. Andruskiewitsch. \emph{An Introduction to Nichols Algebras}. In Quantization, Geometry and Noncommutative Structures in Mathematics and Physics. 
A. Cardona et al., eds., pp. 135--195, Springer-Nature (2017).

\bibitem{aah-triang} N. Andruskiewitsch,  I. Angiono, I. Heckenberger. 
\emph{On finite GK-dimensional Nichols algebras over abelian groups}. Mem. Amer. Math. Soc.
\textbf{271}, No. 1329 (2021).


\bibitem{aah-oddchar} N. Andruskiewitsch,  I. Angiono, I. Heckenberger.
\emph{Examples of finite-dimensional pointed Hopf algebras in positive characteristic}. In
\emph{Representation Theory, Mathematical Physics and Integrable Systems, in honor of N. Reshetikhin}, eds.  A. Alexeev, \emph{et al.} Progr. Math. \textbf{340}, pp. 1--38 (2021).


\bibitem{aapw} N. Andruskiewitsch,  I. Angiono, J. Pevtsova, S. Witherspoon.
\emph{Cohomology rings of finite-dimensional pointed Hopf algebras over abelian groups}.  
Res. Math. Sci. \textbf{9:12} (2022).

\bibitem{andrus-devoto}   N. Andruskiewitsch and J. Devoto. \emph{Extensions of Hopf algebras}.
Algebra i Analiz, \textbf{7} (1)  22--61 (1995); St. Petersburg Math. J. \textbf{7} (1)  17--52 (1995).

\bibitem{ABFF} N. Andruskiewitsch, D. Bagio, S. D. Flora and D. Flores.
\emph{On the Laistrygonian Nichols algebras that are domains}. \texttt{arXiv:2202.13957}.


\bibitem{A-Penha} N. Andruskiewitsch and H. M. Pe\~na Pollastri, 
\emph{On the restricted Jordan plane in odd characteristic}. J. Algebra Appl. \textbf{20} (01), Article No. 2140012 (2021).



\bibitem{A-Schneider-cambr}
N.~Andruskiewitsch and H.-J. Schneider, \emph{Pointed {H}opf algebras}. In \emph{New directions in {H}opf algebras}, pp. 1--68, Math. Sci. Res. Inst. Publ. 43,
Cambridge Univ. Press, Cambridge, (2002). 



\bibitem{bnpp} C.\ Bendel, D.\ K.\ Nakano, B.\ J.\ Parshall and C.\ Pillen. 
\emph{Cohomology for quantum groups via the geometry of the nullcone}.
Mem.\ Amer.\ Math.\ Soc.\ \textbf{229} (2014), no.\ 1077. 


\bibitem{clw} C. Cibils, A. Lauve, S. Witherspoon, \emph{Hopf quivers and Nichols algebras in positive characteristic}, Proc. Amer. Math. Soc. 137(12) (2009) 4029–4041.



\bibitem{drinfeld-quasi} V. G. Drinfeld. \emph{Quasi-Hopf algebras}. Leningrad Math. J. \textbf{1}, pp. 1419--1457 (1990).

\bibitem{drupieski} C. Drupieski.
\emph{Representations and cohomology for Frobenius-Lusztig kernels}. 
J. Pure Appl. Algebra {\bf 215}, pp.  1473--1491 (2011).

\bibitem{drupieski2} C. Drupieski.
\emph{Cohomological finite-generation for finite supergroup schemes}. 
Adv. Math. {\bf 288}, pp. 1360--1432 (2016).

\bibitem{EGNO}
P.~Etingof, S.~Gelaki, D.~Nikshych and V.~Ostrik, \emph{Tensor categories}.
Mathematical Surveys and Monographs 205, Amer. Math. Soc., 2015.


\bibitem{etingof-ostrik03} P. Etingof and V. Ostrik. \emph{Finite
tensor categories}. Mosc. Math. J. {\bf 4}, pp.  627--654, 782--783  (2004).

\bibitem{EOW} K. Erdmann, O. Solberg and X. Wang. 
\emph{On the structure and cohomology ring of connected Hopf algebras}. J. Algebra \textbf{527}, pp. 366--398 (2019).

\bibitem{evens}  L. Evens.  \emph{The cohomology ring of a finite group}.
Trans. Amer. Math. Soc. {\bf 101}, pp. 224--239 (1961).

 

\bibitem{friedlander-negron}  E. Friedlander and C. Negron.  
\emph{Cohomology for Drinfeld doubles of some infinitesimal group schemes}. Alg. Number Th. {\bf 12},  pp. 1281--1309    (2018).


\bibitem{friedlander-parshall}  E. Friedlander and B. Parshall.  
\emph{On the cohomology of algebraic and related finite groups}. Invent. Math. {\bf 74}, pp.  85--117 (1983).

\bibitem{friedlander-suslin} E.\ Friedlander and A.\ Suslin. 
\emph{Cohomology of finite group schemes over a field}.
Invent.\ Math.\ \textbf{127}, pp. 209--270 (1997). 


\bibitem{ginzburg-kumar} V.\ Ginzburg and S.\ Kumar.  \emph{Cohomology
of quantum groups at roots of unity}. Duke Math.\ J.  \textbf{69}, pp.
179--198 (1993).

\bibitem{golod} E. Golod.  \emph{The cohomology ring of a finite $p$-group}.
 Dokl. Akad. Nauk SSSR {\bf 235}, pp. 703--706 (1959).

\bibitem{gordon} I. G. Gordon.  
\emph{Cohomology of quantized function algebras at roots of unity}. 
Proc. London Math. Soc. (3) {\bf 80}, pp. 337–359  (2000). 

\bibitem{hofstetter}   I. Hofstetter,
\emph{Extensions of Hopf algebras and their cohomological description},
J. Algebra {\bf 164} (1994),  246--298.

\bibitem{jacobson} N. Jacobson. \emph{Lie algebras}. 
Reprint of the 1962 original. Dover Publ., Inc., New York, 1979. {\rm ix}+331  


\bibitem{kac} G.~I. Kac, \emph{Extensions of groups to ring groups},
Math. USSR Sbornik {\bf 5} (1968),  451--474.

\bibitem{mackey}
G.~W. Mackey, \emph{Products of subgroups and projective multipliers}, Colloquia Math. Soc.
Janos Bolyai 5, Hilbert Space Operators, Tihany (Hungary) (1970), 401--413.

\bibitem{majid-book}
S.~Majid,  \emph{Foundations of Quantum Group Theory.\ }
Cambridge Univ. Press, 1995.

\bibitem{majid-oeckl}
S.~Majid and R.~Oeckl, \emph{Twisting of quantum differentials and the Planck scale Hopf algebra}, Commun. Math. Phys. \textbf{205}, No. 3, 617--655 (1999).

\bibitem{MPSW} M.\ Mastnak, J. Pevtsova, P. Schauenburg, and S. Witherspoon.
\emph{Cohomology of finite-dimensional pointed Hopf algebras}.
Proc.\ London Math.\ Soc.\ (3) \textbf{100}, pp. 377--404 (2010). 

\bibitem{masuoka-survey} A. Masuoka, \emph{Hopf algebra extensions and cohomology,} In \emph{New directions in {H}opf algebras}, pp. 1--68, 
Math. Sci. Res. Inst. Publ. 43,
Cambridge Univ. Press, Cambridge, (2002). 

\bibitem{muger} M. M\"uger. \emph{From subfactors to categories and topology I. Frobenius algebras in and Morita equivalence of tensor categories}. 
J. Pure Appl. Algebra \textbf{180}, pp. 81--157 (2003).


\bibitem{negron} C. Negron,
\emph{Finite generation of cohomology for Drinfeld doubles of finite group schemes}.
Sel. Math., New Ser. \textbf{27}, No. 2, Paper No. 26, 20 pp. (2021).


\bibitem{negron-plavnik} C. Negron, J. Y. Plavnik.  \emph{Cohomology of finite tensor categories: duality and Drinfeld centers}. Trans. Amer. Math. Soc. \textbf{375} (2022), 2069--2112.

\bibitem{NWW} Van C. Nguyen, X. Wang and S. Witherspoon. \emph{Finite generation of some cohomology rings via twisted tensor product and Anick resolutions}. J. Pure Appl. Algebra \textbf{223}, pp. 316--339 (2019).




\bibitem{radford-book} D. E. Radford.
\emph{Hopf algebras}. Series on Knots and Everything 49. 
Hackensack, NJ: World Scientific. xxii, 559 p.  (2012).

\bibitem{schauenburg-jpaa} P. Schauenburg.
\emph{The monoidal center construction and bimodules}. 
J. Pure Appl. Algebra \textbf{158}, No. 2--3, 325--346 (2001).


\bibitem{schauenburg} P. Schauenburg.
\emph{Hopf bimodules, coquasibialgebras, and an exact sequence of Kac} 
Adv. Math. \textbf{165}, No. 2, 194--263 (2002).



\bibitem{schneider}   H.-J. Schneider,
\emph{A normal basis and transitivity of crossed products for Hopf algebras},
J. Algebra {\bf 152} (1992),  289--312.


\bibitem{stefan-vay}  D. \c{S}tefan and C. Vay.  \emph{The cohomology ring of the 12-dimensional
Fomin–Kirillov algebra}. Adv. Math. \textbf{291}, pp. 584--620 (2016).

\bibitem{sweedler}M.~E. Sweedler, \emph{Hopf algebras},
New York: W.~A. Benjamin, Inc. 1969, 336 p. (1969).


\bibitem{Takeuchi-matched}
M. Takeuchi, \emph{Matched pairs of groups and bismash products of Hopf algebras}, Commun.
Algebra 9 (1981), 841--882.

\bibitem{V} B. B. Venkov.  \emph{Cohomology algebras for some classifying spaces}. Dokl. Akad. Nauk. SSSR {\bf 127},
pp.  943--944  (1959).
\end{thebibliography}
\end{document}